\documentclass[letterpaper,11pt]{article}

\usepackage{ucs}
\usepackage[utf8x]{inputenc}
\usepackage{graphicx}
\usepackage{amsfonts}
\usepackage{dsfont}
\usepackage{amssymb}
\usepackage{amsmath}
\usepackage{amsthm}
\usepackage{enumerate}
\usepackage{stmaryrd}
\usepackage{fullpage}
\usepackage{ifthen}
\usepackage{subfigure}
\usepackage{epic}
\usepackage{authblk}
\usepackage{textcomp}
\usepackage[small]{caption}

\usepackage[hypertexnames=false,colorlinks=true,linkcolor=blue,citecolor=blue]{hyperref}
\usepackage[numbers,comma,square,sort&compress]{natbib}
\usepackage[letterpaper,text={7in,9in},centering]{geometry}

\usepackage{color}
\usepackage{titlesec}
\graphicspath{{eps/}{pdf/}}
\newcommand{\be}{\begin{equation}}
\newcommand{\ee}{\end{equation}}
\newcommand{\bqs}{\begin{equation*}}
\newcommand{\eqs}{\end{equation*}}

\newcommand{\C}{\mathbb{C}}
\newcommand{\N}{\mathbb{N}}
\newcommand{\R}{\mathbb{R}} 
\newcommand{\Z}{\mathbb{Z}}

\newcommand{\A}{\mathcal{A}}
\newcommand{\E}{\mathcal{E}}
\newcommand{\G}{\mathcal{G}}
\newcommand{\K}{\mathcal{K}}
\newcommand{\cR}{\mathcal{R}}

\newcommand{\rme}{\mathbf{e}}

\newcommand{\rmO}{\mathcal{O}}

\newcommand{\e}{\epsilon}

\numberwithin{equation}{section}

\newenvironment{Hypothesis}[1]%
  {\begin{trivlist}\item[]{\bf Hypothesis #1 }\em}{\end{trivlist}}

\theoremstyle{plain}
\newtheorem{theorem}{Theorem}[section]
\newtheorem{proposition}[theorem]{Proposition}
\newtheorem{lemma}[theorem]{Lemma}

\newtheorem{rmk}[theorem]{Remark}

\newcommand{\hv}{\hat{v}}
\newcommand{\hp}{\hat{\phi}}

\newenvironment{Proof}[1][.]%
 {\begin{trivlist}\item[]\textbf{Proof#1 }}%
 {\hspace*{\fill}$\rule{0.3\baselineskip}{0.35\baselineskip}$\end{trivlist}}

\title{Modulated traveling fronts for a nonlocal Fisher-KPP equation: a dynamical systems approach}
\author[1]{Gr\'egory Faye}
\author[2]{Matt Holzer\footnote{Corresponding Author, mholzer@gmu.edu, 703-993-1463}}
\affil[1]{\small CAMS -  Ecole des Hautes Etudes en Sciences Sociales, 190-198 avenue de France, 75013, Paris, France}
\affil[2]{\small Department of Mathematical Sciences, George Mason University, Fairfax, VA 22030, USA}

\begin{document}
\maketitle

\begin{abstract}
\noindent We consider a nonlocal generalization of the Fisher-KPP equation in one spatial dimension.  As a parameter is varied the system undergoes a Turing bifurcation.  We study the dynamics near this Turing bifurcation.  Our results are two-fold.  First, we prove the existence of a two-parameter family of bifurcating stationary periodic solutions and derive a rigorous asymptotic approximation of these solutions. We also study the spectral stability of the bifurcating stationary periodic solutions with respect to almost co-periodic perturbations. Secondly, we restrict to a specific class of exponential kernels for which the nonlocal problem is transformed into a higher order partial differential equation.  In this context, we prove the existence of modulated traveling fronts near the Turing bifurcation that describe the invasion of the Turing unstable homogeneous state by the periodic pattern established in the first part.  Both results rely on a center manifold reduction to a finite dimensional ordinary differential equation.
\end{abstract}


\section{Introduction}

In this paper, we consider the following nonlocal partial differential equation, 
\be
u_t=u_{xx}+\mu u\left(1-\phi\ast u\right), \quad x\in \R, \label{eq:NLKPP}
\ee
where $\mu>0$ represents the strength of the nonlocal competition and
\[ \phi\ast u(x)= \int_\R \phi(x-y)u(y)dy,\]
for a given kernel $\phi$. Throughout this paper, we will assume that the kernel $\phi$ satisfies the following hypotheses.
\begin{Hypothesis}{(H1)} The kernel $\phi$ satisfies:
\[   \phi\geq 0, \quad \phi(0)>0, \quad  \phi(-x)=\phi(x), \quad \int_\R \phi(x)dx=1, \text{ and } \int_\R x^2\phi(x)dx<\infty.\]
\end{Hypothesis}
In the limiting case where $\phi$ is replaced by the Dirac $\delta$-function, the nonlocal partial differential equation \eqref{eq:NLKPP} reduces to the classical Fisher-KPP equation \cite{Fis,Kol-Pet-Pis}
\begin{equation}\label{kpp-local}
u_t=u_{xx}+\mu u(1-u), \quad x\in\R.
\end{equation}
Such an equation \eqref{kpp-local} arises naturally in many mathematical models in biology, ecology or genetics, see \cite{Fis,Kol-Pet-Pis}, and $u$ typically stands the density of some population. The nonlocal equation \eqref{eq:NLKPP} can then be interpreted as a generalization of the local Fisher-KPP equation \eqref{kpp-local} in which interactions among individuals are nonlocal. For more details on such nonlocal models, we refer to \cite{Fur-Gri,Gou,Gen-Vol-Aug} among others.

\paragraph{History of the problem.}

The behavior of the solutions of the local equation \eqref{kpp-local} has been studied for decades and is now well understood while much less is known about solutions to the nonlocal equation \eqref{eq:NLKPP}. Indeed, from a mathematical point of view, the analysis of \eqref{eq:NLKPP}  is quite involved since this class of equations with a nonlocal competition term generally does not satisfy the comparison principle. Recently, theoretical and numerical studies \cite{Ber-Nad-Per-Ryz,Fan-Zha, Nad-Per-Tan} have shown that for sufficiently small $\mu$, the solutions share many of the same properties of the local Fisher-KPP equation in that there exists a family of traveling wave solutions of the form
\be u(t,x)=U(x-ct), \quad \underset{\xi \rightarrow -\infty}{\lim}U(\xi)=1, \quad \underset{\xi \rightarrow +\infty}{\lim}U(\xi)=0, \quad U \text{ decreasing}.\label{eqFront} \ee
It is known that, again for  $\mu$ sufficiently small, these traveling waves and the homogeneous stationary solutions $u(t,x)=1$ and $u(t,x)=0$ are the only bounded solutions to \eqref{eq:NLKPP}, see \cite{Ach-Kue,Alf-Cov,Ber-Nad-Per-Ryz,Fan-Zha}. On the other hand, when $\mu$ is large, some other bounded solutions may exist as suggested by the numerical exploration of \cite{Nad-Per-Tan}. More precisely, if the Fourier transform of the kernel $\phi$ takes some negative values, then for sufficiently large $\mu$, the trivial state $u(t,x)=1$ is Turing unstable for \eqref{eq:NLKPP}.  This suggests the emergence of non-monotonic bounded solutions \cite{Apr-Bes-Vol-Vou,Gen-Vol-Aug}.  Indeed, recent work by Hamel and Ryzhik \cite{Ham-Ryz} has shown the existence of stationary periodic solutions $u$ satisfying,
\be 0=u_{xx}+\mu u(1-\phi*u), \quad x\in \R \label{eqPer}\ee
for large $\mu$ when the Fourier transform of the kernel attains negative values.

\begin{figure}[ht]
\centering
   \includegraphics[width=0.4\textwidth]{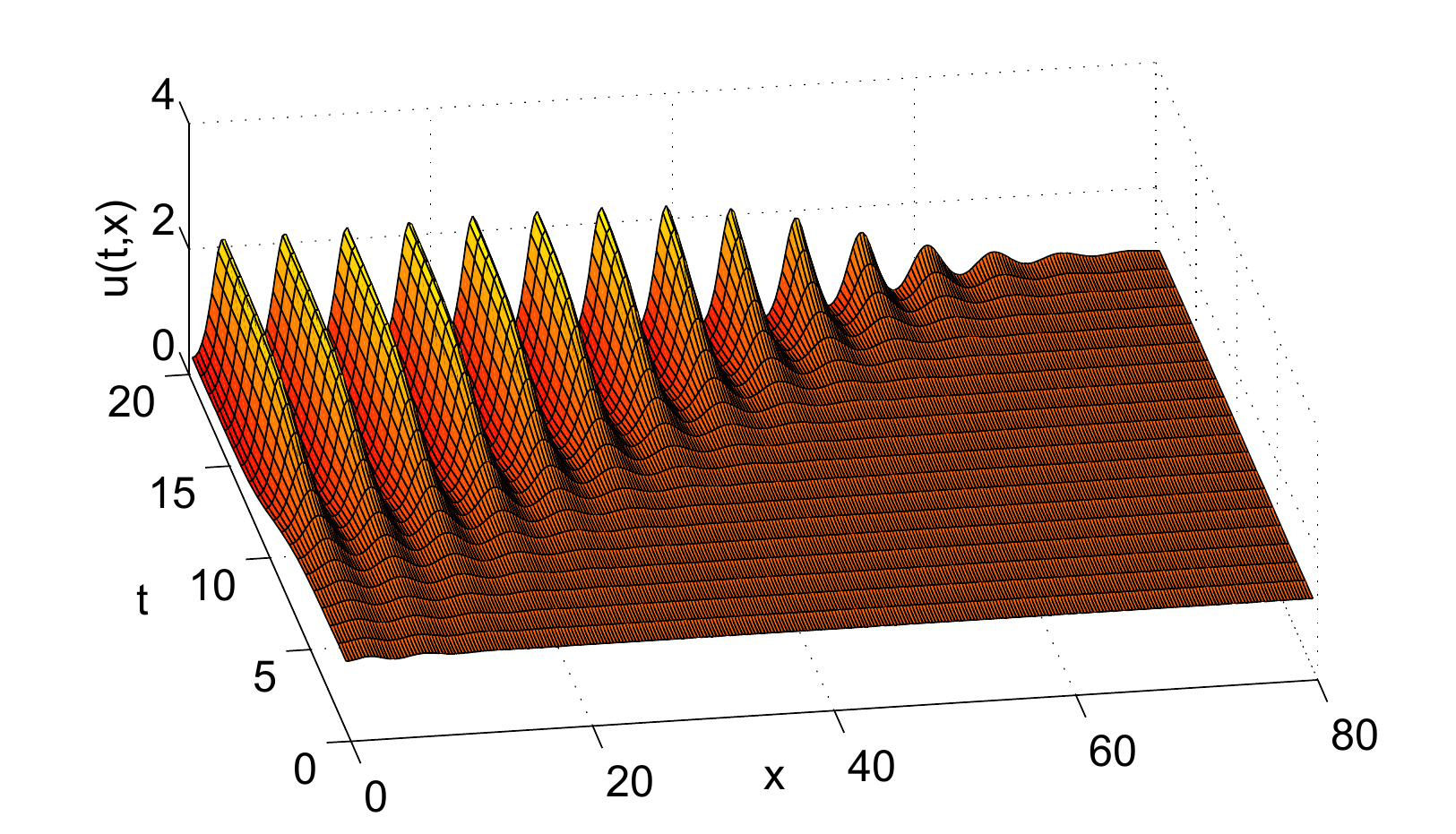}
   \includegraphics[width=0.4\textwidth]{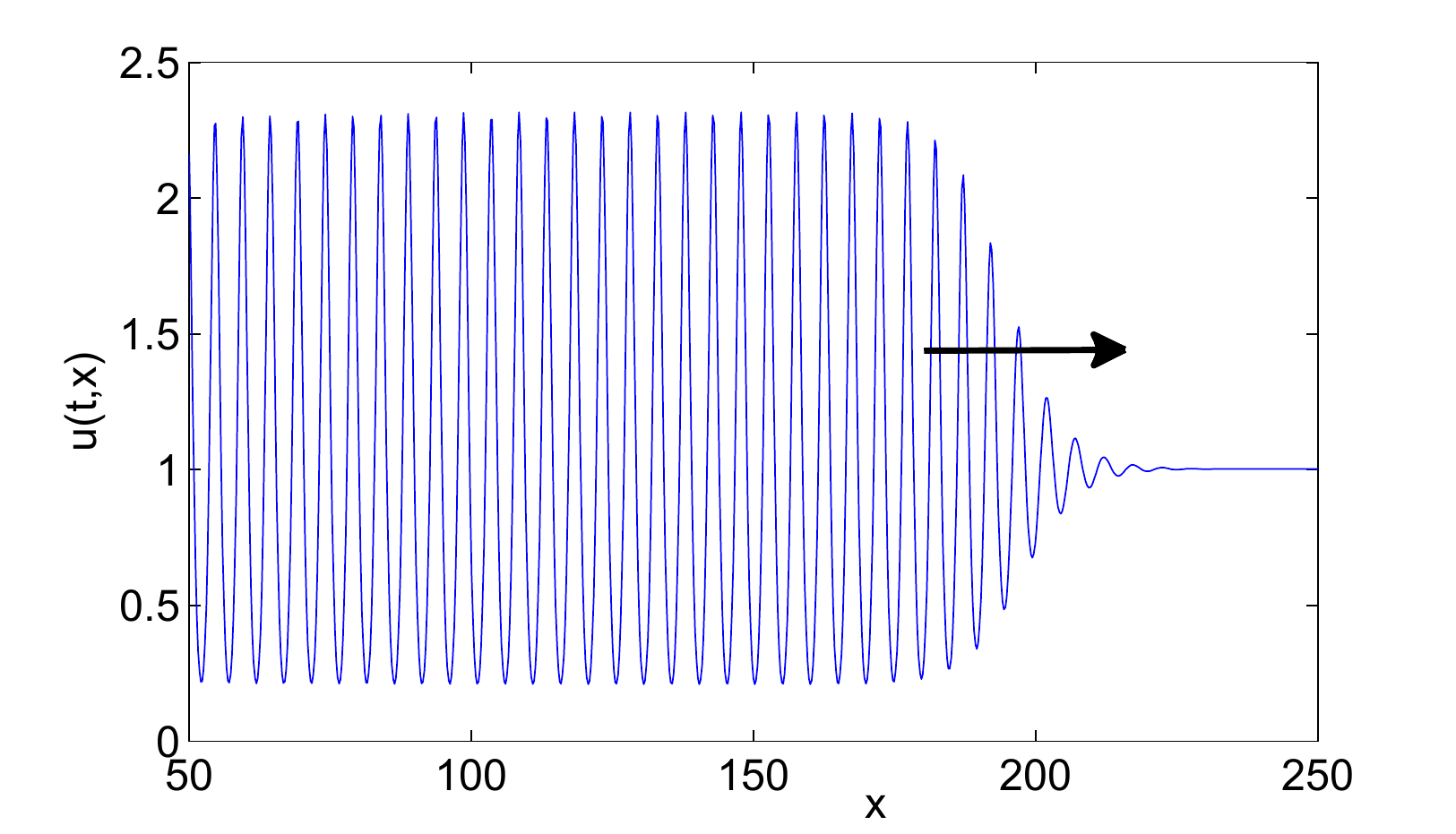}  
\caption{A modulated traveling front obtained from direct numerical simulation of  (\ref{eq:NLKPP}) with kernel (\ref{eq:kernel}) for $a=0.7$ and $\mu=32$.  The front invades the Turing unstable state $u=1$ and leaves a stationary periodic pattern in its wake.     }
\label{fig:spacetime}
\end{figure}

Recent numerical studies of (\ref{eq:NLKPP}) also suggest the existence of modulated traveling fronts where these stationary periodic solutions invade the Turing unstable state $u=1$, see \cite{Nad-Per-Tan} and Figure~\ref{fig:spacetime}.  These modulated traveling fronts are the focus of study of this article.  For a certain class of kernels, we will prove the existence of modulated traveling fronts of the form
\be u(t,x)=U(x-ct,x), \quad \underset{\xi \rightarrow -\infty}{\lim}U(\xi,x)=1+P(x), \quad \underset{\xi \rightarrow +\infty}{\lim}U(\xi,x)=1,\label{eqModTF}\ee 
where $P(x)$ is a stationary periodic solution of the shifted problem
\be 0=v_{xx}-\mu v -\mu v\phi*v, \quad x\in \R \label{eqPerMod}.\ee
The first step of our analysis will be to refine the existence result of \cite{Ham-Ryz} for parameter values near the onset of Turing instability and then use center manifold techniques to construct modulated traveling fronts of the form \eqref{eqModTF}. In that direction, we also point out that explicit examples of wave-train solutions have been recently constructed in \cite{Duc-Nad,Nad-Per-Ros-Ryz} for a different nonlocal problem.

\paragraph{The assumptions.}
Before stating our main results, we first make some further assumptions on the kernel $\phi$. 
Linearizing equation \eqref{eq:NLKPP} around the stationary homogeneous state $u=1$, we find the following dispersion relation,
\be d(\lambda,k,\mu):=-k^2-\mu \hat \phi(k)-\lambda.
\label{eq:disp}\ee

\begin{Hypothesis}{(H2)}
For $\phi(x)$ satisfying (H1),  we further assume that there exists unique $k_c>0$ and $\mu_c>0$, such that the following conditions are satisfied,
\begin{itemize}
\item[(i)]  $d(0,k_c,\mu_c)=0$.  
\item[(ii)] $\partial_k d (0,k_c,\mu_c)=0$.
\item[(iii)] $\partial_{kk} d (0,k_c,\mu_c)<0$.
\end{itemize}
\end{Hypothesis}
The first condition imposes that $\hat \phi(k_c)<0$ as from the dispersion relation $d(0,k_c,\mu_c)=0$, we have that $\hat \phi(k_c)=-\frac{k_c^2}{\mu_c}<0$. The second condition ensures that $k_c$ is a double root of the dispersion relation and combined with the third condition, that $\mu_c$ represents the onset of instability.  For $\mu<\mu_c$ all the spectrum is to the left of the imaginary axis while for $\mu>\mu_c$ there is a band of wavenumbers surrounding $k=\pm k_c$ that are unstable.  

For some part of our analysis, we will work with a specific kernel that satisfies all the hypotheses (H1), namely we will choose
\be 
\label{eq:kernel}
\phi(x):=Ae^{-a|x|}-e^{-|x|},
\ee
for some values of $A>0$ and $a>0$.  Recall that (H1) requires that $\phi(x)>0$ and $\int_{-\infty}^\infty \phi(x)dx=1$.  The second condition implies that $A=3a/2$, and the first condition in turn implies that $a\in\left(2/3,1\right)$.  The choice of such a specific kernel is motivated by the fact that equation \eqref{eq:NLKPP} can be reduced to a system of partial differential equations. Indeed, define
\[ v(t,x):=Ae^{-a|x|}\ast u(t,x), \quad w(t,x):=-e^{-|x|}\ast u(t,x),\]
and let $\phi_v:=Ae^{-a|x|}$ and $\phi_w:=-e^{-|x|}$ for future reference. We find that (\ref{eq:NLKPP}) reduces to the following system,
\begin{subequations}
\begin{align}
u_t &= u_{xx}+\mu u (1-v-w),  \\
0 &= v_{xx}-a^2 v+3a^2u, \\
0&= w_{xx}-w-2u.
\end{align}
\label{systemPDE}
\end{subequations}
Using a specific connectivity kernel can sound very restrictive, but it has proven in other contexts its efficiency to overcome the difficulty of the nonlocal nature of the problem while still gaining some general insights. In particular, we refer to some recent works on the existence and stability of traveling pulses in neural field equations with synaptic depression or on some pinning and unpinning phenomena in nonlocal systems \cite{Pinning,Fay} where kernels with rational Fourier transform have been used to reduced the problem to a high-order system of partial differential equations.
 
%\caption{Numerical simulation of (\ref{eq:NLKPP}) with kernel (\ref{eq:kernel}).  Here we have used $a=0.7$ and $\mu=32$. Initial conditions are identically one with a compactly supported perturbation placed at the left boundary of the domain.  On the left is a space-time plot of the solution for $t\in[0,20]$.  On the right, is the solution at fixed time $t=48$.  Observe the traveling front propagating to the right that leaves a stationary periodic pattern in its wake.   }

\paragraph{Main results.} This paper contains three main results. The first result concerns the existence of stationary periodic solutions of the nonlocal equation \eqref{eq:NLKPP} and can be stated as follows.

\begin{theorem}\label{thmsolper} Assume that hypotheses $(H1)$ and $(H2)$ are satisfied.  Let $\mu:=\mu_c+\e^2$ and $k:=k_c+\delta$. There exists $\epsilon_0>0$, such that for all $\e \in(0,\e_0]$ and all $\delta^2 < \dfrac{-\hat \phi(k_c)}{1+\frac{\mu_c}{2}\hat\phi''(k_c)}\e^2$ there is a stationary $\frac{2\pi}{k}$-periodic solution of \eqref{eq:NLKPP} with leading expansion of the form
\be \mathbf{u}_{\e,\delta}(x) = 1+ \sqrt{\frac{\hat\phi(k_c)\epsilon^2+\left(1+\frac{\mu_c}{2}\hat\phi''(k_c)\right)\delta^2}{\omega}} \cos\left((k_c+\delta)x\right) + \rmO\left( \left| \epsilon^2-\delta^2 \right|\right), \label{solper}\ee
where $\omega<0$ is defined in equation \eqref{eqcoeffomega}. Moreover, for any $\tau \in [0,2\pi/k]$, $(x)\mapsto \mathbf{u}_{\e,\delta}(x+\tau)$ is also a solution of \eqref{eq:NLKPP}.
\end{theorem}

First, note that the results of Theorem \ref{thmsolper} do not rely on a specific form of the kernel. This theorem complements the study of Hamel \& Ryzhik \cite{Ham-Ryz} where they also proved the existence of stationary periodic solutions of \eqref{eq:NLKPP}. While the analysis in \cite{Ham-Ryz} is global and relies on degree theory and in the regime $\mu$ large, our study is local and uses center manifold theory. To some extent, our approach gives sharper results close to the bifurcation point $\mu_c$ as we obtain a complete description of all bounded stationary solutions of \eqref{eq:NLKPP} in some neighborhood of the solution $u=1$. Furthermore, we show the existence of a family of periodic solutions indexed by their spatial frequency $k\approx k_c$. 

Our second main result is a spectral analysis of the stationary periodic solutions found in Theorem \ref{thmsolper} and our results are summarized as follows.

\begin{theorem}\label{thmstability}
Assume that hypotheses $(H1)$ and $(H2)$ are satisfied. Then, the following assertions are true.
\begin{itemize}
\item[(i)] The periodic solutions $\mathbf{u}_{\e,\delta}$ given in Theorem \ref{thmsolper} are neutrally stable with respect to perturbations of the same period $2\pi/k$, where $k=k_c+\delta$ and $\delta$ satisfies the relation $\delta^2 < \dfrac{-\hat \phi(k_c)}{1+\frac{\mu_c}{2}\hat\phi''(k_c)}\e^2$ for all $\e \in(0,\e_0]$.
\item[(ii)] The periodic solutions $\mathbf{u}_{\e,\delta}$ given in Theorem \ref{thmsolper} are spectrally unstable with respect to perturbations of the form $e^{i\sigma x}V(x)$, where $V\in \mathcal{Y}$ is a solution of the spectral problem \eqref{sepctralPb}, in the limit $\sigma \rightarrow 0$ and whenever $\delta$ satisfies
\be \frac{-\hat \phi(k_c)}{3\left(1+\frac{\mu_c}{2}\hat\phi''(k_c)\right)} \e^2 < \delta^2< \frac{-\hat \phi(k_c)}{1+\frac{\mu_c}{2}\hat\phi''(k_c)}\e^2 ,\label{eqStabRegion}
\ee
for $\e\in(0,\e_0]$.
\end{itemize}
\end{theorem}

The first part of this theorem is a direct consequence of the center manifold reduction used in the existence proof of Theorem \ref{thmsolper}. Indeed, the spectral analysis of the periodic solutions with respect to perturbations of the same period can be directly done on the reduced two dimensional  equation on the center manifold where one finds two eigenvalues $\lambda_c=0$ and $\lambda_s<0$. The fact that there exists a critical eigenvalue is due to the translation of invariance of the problem, namely $\frac{\mathrm{d}}{\mathrm{d}x}\mathbf{u}_{\e,\delta}$ is always in the kernel of the linearized operator. The second part of the theorem is a perturbation analysis, using Liapunov Schmidt reduction, where we show that the critical eigenvalue $\lambda_c=0$ is perturbed into $\lambda_c = \mathbf{g}(\e,\delta) \sigma^2+\rmO(\sigma^4)$, as $\sigma \rightarrow 0$ when we are looking for perturbations of the form $e^{i\sigma x}V(x)$, $V\in \mathcal{Y}$. The region in parameter space $(\e,\delta)$ where $\mathbf{g}(\e,\delta)>0$ will then give spectral instability with respect to such perturbations. Finally, note that instabilities with respect to perturbations with almost same frequencies are called {\it sideband instabilities} in the literature \cite{Mie1}.

The central result of this paper is the proof of the existence of modulated traveling front solutions that are asymptotic at infinity to the stationary periodic solutions found in Theorem \ref{thmsolper} and the homogeneous state $u=1$. We will realize these modulated traveling fronts as heteroclinic orbits of a reduced system of ODEs in normal form. Roughly speaking, and fixing the frequency to $k_c$, we look for solutions of \eqref{eq:NLKPP} that can be written as
\[ u(t,x) = U(x-\e st,x)=\sum_{n\in \Z} U_n(x-\e st)e^{-in k_cx}, \]
where $\e=\sqrt{\mu-\mu_c}$. Replacing this ansatz into the equivalent system \eqref{systemPDE} will lead to the study of an infinite dimensional dynamical system (see Section \ref{secModTF} and equation \eqref{eqTWmode} for more details) of the form
\be
\label{eqIntro}
\partial_\xi U_n = \mathcal{L}_n^{\e} U_n + \mathcal{R}_n(U,\e).
\ee
The main difficulty in studying \eqref{eqIntro}  comes from the presence at onset ($\e=0$) of an infinite dimensional central part. However, as $0<\e\ll1$, these eigenvalues will leave the imaginary axis with different velocities. A finite number will stay close $(\rmO(\e))$ to the imaginary axis while all other eigenvalues leave fast enough ($\rmO(\sqrt{\e})$) so that a spectral gap exists.  This gap will allow for small $\e>0$ the construction of a finite dimensional invariant manifold of size $\rmO(\e^{2/3+\gamma})$, for $\gamma>0$. This manifold will contain the modulated traveling fronts that we are looking for. The result that we obtain can be formulated as follows.

\begin{theorem}\label{thmModTF}
Assume that $\phi$ is the kernel given in \eqref{eq:kernel} and that hypothesis (H2) is satisfied. Provided that $s>\sqrt{-4\hat\phi(k_c)\zeta}$, where $\zeta$ is defined in \eqref{eqcoeffzeta}, there is an $\e_0>0$ such that for all $\e\in (0,\e_0)$, and all $\delta^2< -\frac{\hat\phi(k_c)}{1+\frac{\mu_c}{2}\hat\phi''(k_c)}\e^2$, equation \eqref{eq:NLKPP} has modulated traveling front solutions of frequency $k_c+\delta$ and of the form 
\[ u(t,x) = U(x-\e st,x)=\sum_{n\in \Z} U_n(x-\e st)e^{-in(k_c+\delta)x}, \]
with the boundary conditions at infinity
\[ \underset{\xi\rightarrow - \infty}{\lim}U(\xi,x) = \mathbf{u}_{\e,\delta}(x) \text{ and } \underset{\xi\rightarrow + \infty}{\lim}U(\xi,x) =1.\]
\end{theorem}

 The first known existence results of modulated traveling waves are due to Collet \& Eckmann \cite{Col-Eck1,Col-Eck2} and Eckmann \& Wayne \cite{Eck-Way}, who proved  the existence of such solutions in the  Swift-Hohenberg equation with cubic nonlinearities. The techniques developed in \cite{Eck-Way} have then been generalized for the problem of bifurcating fronts for the Taylor-Couette problem in infinite cylinders by Haragus \& Schneider in \cite{Har-Sch} with quadratic nonlinearities, and our proof of Theorem \ref{thmModTF} will rely on a center manifold result presented in \cite{Har-Sch}. Finally, note that similar results have been obtained in the two-dimensional Swift-Hohenberg equation for more general modulated fronts, for example modulated fronts that connect stable hexagons with unstable roll solutions \cite{Doe-San-Sch-Sch}.

\paragraph{Outline of the paper.} The paper is organized as follows. In section \ref{secSatPer}, we prove the existence of a two-parameter family of stationary periodic solutions and consider their stability with respect to almost co-periodic perturbations. Then in section \ref{secAmpEq}, we describe the modulated fronts for which we are looking for  as solutions of an associated amplitude equation. Finally, in the last section \ref{secModTF}, we prove the existence of modulated traveling fronts.  We conclude in section~\ref{sec:discussion} with an discussion of future directions for study.

\section{Stationary Periodic Solutions}\label{secSatPer}

In this section, we first prove the existence of stationary periodic solutions of equation \eqref{eq:NLKPP} for $\mu \approx \mu_c$ and $k\approx k_c$. In a second step, we study their instabilities with respect to perturbations with different but close frequencies.

\subsection{Existence of stationary periodic solutions}

Throughout the paper, we use the notations $\mu:=\mu_c+\e^2$ and $k:=k_c+\delta$, with $0<\e\ll 1$ and $0<\delta \ll 1$. In order to prove the existence of stationary periodic solutions of equation \eqref{eq:NLKPP}, we first rescale space and change coordinates, $u=1+v$, such that we obtain the new equation
\[
v_t=k^2 v_{xx}-\mu \phi_k*v -\mu v\phi_k * v,\quad \phi_k (x) =\frac{1}{k}\phi\left(\frac{x}{k}\right) ,
\]
where $v$ is $2\pi$-periodic in $x$. This equation can formally be written as
\be
v_t = \A(\mu,k) v+\cR(v,\mu,k), \label{eqFormal}
\ee
where
\[ \A(\mu,k) v:= k^2 v_{xx}-\mu \phi_k*v, \quad \cR(v,\mu,k):=-\mu v\phi_k * v.\]
If we define,
\begin{align*} \mathcal{X}&:=L^2_{per}[0,2\pi]=\left\{ u \in L^2_{loc}(\R)~|~ u(x+2\pi)=u(x), \quad x\in\R\right\}, \\
 \mathcal{Y}&:=\mathcal{D}(\A)=H^2_{per}[0,2\pi]=\left\{ u \in H^2_{loc}(\R)~|~ u(x+2\pi)=u(x), \quad x\in\R\right\}, \end{align*}
then we readily see that the nonlinear map $\cR: \mathcal{Y} \rightarrow \mathcal{Y}$ is smooth. We now rewrite \eqref{eqFormal} as
\be
v_t = \A_c v+\mathcal{B}(\epsilon,\delta)v+ \cR(v,\mu_c+\e^2,k_c+\delta), \label{eqFormalMod}
\ee
with 
\[ \A_c:=\A(\mu_c,k_c)\text{ and }  \mathcal{B}(\epsilon,\delta):= \A(\mu_c+\e^2,k_c+\delta)-\A_c.\]
We note that $\A_c:\mathcal{Y}\rightarrow\mathcal{X}$ is a continuous linear operator and $\mathcal{Y}$ is dense in and compactly embedded into $\mathcal{X}$, we conclude that the resolvent of $\A_c$ is compact and thus its spectrum $\sigma(\A_c)$ only consists of eigenvalues $\lambda$. According to the dispersion relation \eqref{eq:disp} we have
%\[ \lambda_\ell = -k_c^2\ell^2-\mu_c  \left(\frac{3a^2}{a^2+k_c^2\ell^2}-\frac{2}{1+k_c^2\ell^2}\right), \quad \ell \in \Z, \]
\be\sigma(\A_c)=\left\{ \lambda_\ell \in \C~|~\lambda_\ell = -k_c^2\ell^2-\mu_c \hat{\phi}(k_c\ell), \quad \ell \in \Z\right\}.\label{eqSpectrumAc} \ee
As a consequence,
\[ \sigma(\A_c)\cap i\R=\left\{ 0\right\}, \]
and $\lambda=0$ is an eigenvalue with geometric multiplicity two and associated eigenvectors $\rme(x):=e^{ix}$ and $\bar\rme(x)=e^{-ix}$. One can check that the algebraic multiplicity is also two. We define $\mathcal{X}_c:=\left\{ \rme,\bar\rme\right\}$ and the spectral projection $\mathcal{P}_c:\mathcal{X}\rightarrow \mathcal{X}_c$  via
\[\mathcal{P}_cu=\langle u, \rme \rangle \rme + \langle u, \bar\rme \rangle \bar\rme, \]
where 
\[ \langle u,v \rangle =\frac{1}{2\pi} \int_0^{2\pi} u(x) \bar v(x) dx.\]
From \eqref{eqSpectrumAc}, we see that the spectrum of $\A_c |_{(\text{id}-\mathcal{P}_c)\mathcal{X}}$ is off the imaginary axis, and one can then deduce that
\[ \left\| (i\nu - \A_c)^{-1} \right\|_{(\text{id}-\mathcal{P}_c)\mathcal{X}}\leq \frac{C}{1+|\nu|}, \quad \nu\in\R, \]
with some positive constant $C>0$. Therefore, by applying the center manifold theorem (see \cite{Har-Ioo}, \S 2), we have the existence of neighborhoods of the origin $\mathcal{U}\subset \mathcal{X}_c$, $\mathcal{V}\subset (\text{id}-\mathcal{P}_c)\mathcal{Y}$, $\mathcal{W}\subset \R^2$ and, for any $m<\infty$, a $\mathcal{C}^m$-map $\Psi:\mathcal{U}\times \mathcal{W}\rightarrow \mathcal{V}$ having the following properties.
\begin{itemize}
\item[(i)] For any $(\epsilon,\delta)\in\mathcal{W}$, all bounded solutions $v$ of \eqref{eqFormal} within $\mathcal{U}\times \mathcal{V}$ are on the center manifold, i.e., 
\be
v(t) = A(t) \rme + \bar A(t) \bar \rme + \Psi(A(t),\bar A(t), \epsilon, \delta),\quad \forall t \in \R.
\ee
\item[(ii)] The center manifold is tangent to the center eigenspace, i.e., 
\[ \left\| \Psi(A,\bar A, \epsilon, \delta) \right\| _\mathcal{Y}=\rmO\left( |\epsilon|^2|A|+|\delta| |A|+|A|^2 \right). \]
\item[(iii)] The action of the translation symmetry $\left(\mathbf{T}_\tau u\right)(x)=u(x+\tau)$ and the reflection symmetry $\left(\mathbf{R} u\right)(x)=u(-x)$ on \eqref{eqFormal} implies that the reduced vector field on the center manifold can be written as
\be \frac{dA}{dt}=f(A,\bar A, \epsilon, \delta)=A g(|A|^2,\e,\delta), \label{eqCM} \ee
where $g$ is a $\mathcal{C}^{m-1}$-map in $(A,\bar A,\e,\delta)$ and is real-valued. \end{itemize} 

\begin{lemma}
The Taylor expansion of the map $g$ is given by
\be g(|A|^2,\e,\delta) =  -\hat\phi(k_c)\epsilon^2-\left(1+\frac{\mu_c}{2}\hat\phi''(k_c)\right)\delta^2+\omega |A|^2 +\rmO\left( |\delta|^3+|\epsilon||\delta|+|A|^4 \right) \label{Taylorg} \ee
where
\be \omega:= \mu_c\hat{\phi}(k_c)\left(\frac{\mu_c(\hat{\phi}(k_c)+\hat{\phi}(2k_c))}{4k_c^2+\mu_c\hat{\phi}(2k_c)}+2(1+\hat{\phi}(k_c))\right)<0. \label{eqcoeffomega}\ee
\end{lemma}

\begin{Proof}
Substituting $v(t)=A(t) \rme + \bar A(t) \bar \rme + \Psi(A(t),\bar A(t), \epsilon, \delta)$ into the right-hand side of \eqref{eqFormal} and collecting only the linear terms in $A$, one obtains
\[\frac{dA}{dt}=\left(-(k_c+\delta)^2-(\mu_c+\epsilon^2)\hat \phi(k_c+\delta)\right)A+\rmO_{(\epsilon,\delta)}(A|A|^2). \]
The criticality of the dispersion relation $d(\lambda,k,\mu)=0$ (see \eqref{eq:disp}) implies that
\begin{align*}
-k_c^2-\mu_c\hat \phi(k_c)&=0,\\
-2k_c-\mu_c\hat\phi'(k_c)&=0.
\end{align*}
As a consequence, we have
\begin{align*}
\lambda(\e,\delta)&:=-(k_c+\delta)^2-(\mu_c+\epsilon^2)\hat \phi(k_c+\delta)\\
&=-\left((k_c+\delta)^2-k_c^2-2k_c\delta\right)-\mu_c\left(\hat \phi(k_c+\delta)-\hat\phi(k_c)-\delta\hat\phi'(k_c)\right)-\e^2\hat\phi(k_c+\delta)\\
&=\left(-1-\frac{\mu_c}{2}\hat\phi''(k_c)\right)\delta^2-\hat \phi(k_c)\e^2+\rmO\left(|\e|^2|\delta|+|\delta|^3 \right),
\end{align*}
as $(\e,\delta)\rightarrow(0,0)$.

In order to compute the coefficient $\omega$ in \eqref{Taylorg}, we set $(\e,\delta)=(0,0)$ into equation \eqref{eqFormal} and look for solutions that can be Taylor expanded as
\[v=A \rme + \bar A \bar \rme + A^2\rme_{2,0}+A\bar A\rme_{1,1}+\bar A^2\rme_{0,2}+\rmO(|A|^3). \]
Upon replacing these solutions into \eqref{eqFormal}, we find a hierarchy of equations in powers of $A$ and $\bar A$ that can be solved. The coefficient $\omega$ is readily obtained by projection with $\mathcal{P}_c$ and one finds that
\[ \omega = 2 \langle \widetilde{\cR}(\bar \rme, \rme_{2,0}) + \widetilde{\cR}( \rme, \rme_{1,1}) ,\rme \rangle, \]
where $\widetilde{\cR}$ is a bilinear map defined as
\[ \widetilde{\cR}(u, v):= -\frac{\mu_c}{2}\left( u \phi * v + v \phi * u\right), \quad \forall (u,v)\in \mathcal{Y}^2.  \]
Straightforward computations show that
\begin{align*}
\rme_{2,0}(x)&=\frac{\mu_c\hat \phi(k_c)}{-4k_c^2-\mu_c\hat\phi(2k_c)}e^{i2x}+\text{Span}\left(\rme,\bar\rme \right),\\
\rme_{1,1}(x)&=-2\frac{\hat \phi(k_c)}{\hat \phi(0)}+\text{Span}\left(\rme,\bar\rme \right).
\end{align*}
Replacing these values into the definition of $\omega$ concludes the proof.
\end{Proof}

\begin{Proof}{[of Theorem \ref{thmsolper}]} First we define $\Gamma>0$ as
\[ \Gamma := \frac{\hat\phi(k_c)\epsilon^2+\left(1+\frac{\mu_c}{2}\hat\phi''(k_c)\right)\delta^2}{\omega}.  \]
We are interested in nontrivial stationary solutions $A_0\in \C$ of \eqref{eqCM} that are solutions of 
\[ 0 = g(|A|^2,\e,\delta). \]

Rescaling $A_0 = \sqrt{\Gamma}\widetilde A_0$ and substituting into this equation, we obtain
\[ 0 =  \Gamma\left(  -\omega + \omega |\widetilde A_0|^2 + \rmO(\sqrt{\Gamma}) \right), \text{ as } \Gamma\rightarrow 0. \]
Using the implicit function theorem, one finds solutions of the form 
\[ |\widetilde{A}_0|=1+\rmO(\sqrt{\Gamma}), \text{ as } \Gamma\rightarrow 0. \]
As a conclusion, one has the existence of periodic solutions of \eqref{eq:NLKPP} that can be written as
\[ \mathbf{u}_{\e,\delta}(x) =1+ \sqrt{\frac{\hat\phi(k_c)\epsilon^2+\left(1+\frac{\mu_c}{2}\hat\phi''(k_c)\right)\delta^2}{\omega}} \cos\left((k_c+\delta)x\right) + \rmO\left( \left| \hat\phi(k_c)\epsilon^2+\left(1+\frac{\mu_c}{2}\hat\phi''(k_c)\right)\delta^2 \right|\right) \]
for small $\e \in(0,\e_0]$ and $\delta^2 < \dfrac{-\hat \phi(k_c)}{1+\frac{\mu_c}{2}\hat\phi''(k_c)}\e^2$.
\end{Proof}

\subsection{Sideband instabilities}

This section is devoted to the proof Theorem \ref{thmstability}. In particular we will show the existence of regions in parameter space $(\e,\delta)$ where the periodic solutions found in the previous section are spectrally unstable. The proof involves a Liapunov Schmidt reduction relying on the spectral projection $\mathcal{P}_c$ defined in the previous section. Such a technique was introduced by Mielke for the study of sideband instabilities in the Swift-Hohenberg equation \cite{Mie1}.

First, we denote by $\mathbf{v}_{\e,\delta}$ the periodic solution of \eqref{eqFormal} with expansion
\[\mathbf{v}_{\e,\delta}(x) = \sqrt{\frac{\hat\phi(k_c)\epsilon^2+\left(1+\frac{\mu_c}{2}\hat\phi''(k_c)\right)\delta^2}{\omega}} \cos\left((k_c+\delta)x\right) + \rmO\left( \left| \hat\phi(k_c)\epsilon^2+\left(1+\frac{\mu_c}{2}\hat\phi''(k_c)\right)\delta^2 \right|\right)\]
for small $\e \in(0,\e_0]$ and $\delta^2 < \dfrac{-\hat \phi(k_c)}{1+\frac{\mu_c}{2}\hat\phi''(k_c)}\e^2$. We then linearize \eqref{eqFormal} around $\mathbf{v}_{\e,\delta}$ to obtain the equation
\be
\label{eqLinearized}
v_t = \mathcal{A}(\mu,k)v+D_v\mathcal{R}(\mathbf{v}_{\e,\delta},\mu,\delta)v.
\ee
To show spectral instability of $\mathbf{v}_{\e,\delta}$, we will show that 
\[
\lambda v = \mathcal{A}(\mu,k)v+D_v\mathcal{R}(\mathbf{v}_{\e,\delta},\mu,\delta)v
\]
has a solution $(\lambda,v)$ with $\Re(\lambda)>0$ and $v\neq0$. Here, we allow $v$ to be in $W^{2,\infty}(\R)$ rather than $\mathcal{Y}$. More precisely, we look for solutions of \eqref{eqLinearized} of the form $v(t,x)=e^{\lambda t + i\sigma x}V(x)$, with $V\in \mathcal{Y}$. Thus we arrive at the spectral problem
\be
0 = \G(\e,\delta,\sigma,\lambda) V:= k^2\left(\partial_x+i\sigma \right)^2V-\mu \K(k,\sigma)\cdot V-\mu\left( \mathbf{v}_{\e,\delta} \K(k,\sigma)\cdot V + \phi_{k}*\mathbf{v}_{\e,\delta}V  \right) -\lambda V,
\label{sepctralPb}
\ee
where $k=k_c+\delta$, $\mu=\mu_c+\e^2$ and 
\[ \K(k,\sigma)\cdot V(x) := e^{-i\sigma x} \phi_k *\left(e^{i\sigma \cdot }V \right)(x)=\frac{1}{2\pi}\int_\R \hat \phi (k(\ell+\sigma)) \hat V(\ell)e^{i\ell x}d\ell, \quad x\in\R. \]
We readily note that $\G(0,0,0,0)$ coincides with the linear operator $\mathcal{A}_c$. As a consequence, for sufficiently small $(\e,\delta,\sigma,\lambda)$, we will solve $0 = \G(\e,\delta,\sigma,\lambda) V$ using a Liapunov Schmidt reduction with the splitting we used for the center manifold reduction, i.e. $V$ is decomposed as
\[ V = A \rme +\bar A \bar \rme + \widetilde{V}, \quad \mathcal{P}_c \widetilde{V}=0. \]
Then, for sufficiently small $(\e,\delta,\sigma,\lambda)$, the equation $\left(\text{id}-\mathcal{P}_c\right) \G(\e,\delta,\sigma,\lambda) V=0$ can be solved uniquely using the implicit function theorem for $\widetilde{V}=A v(\e,\delta,\sigma,\lambda)+\bar A \bar v(\e,\delta,\sigma,\lambda)$, with 
\[|v(\e,\delta,\sigma,\lambda)|=\rmO\left(\left| \hat\phi(k_c)\epsilon^2+\left(1+\frac{\mu_c}{2}\hat\phi''(k_c)\right)\delta^2 \right|^{1/2} \right).\]
We can now replace this expression for $\widetilde{V}$ into $\mathcal{P}_c \G(\e,\delta,\sigma,\lambda) V=0$ to get a reduced eigenvalue problem of the form
\[
\mathcal{P}_c \G(\e,\delta,\sigma,\lambda) \left( A \rme +\bar A \bar \rme + A v(\e,\delta,\sigma,\lambda)+\bar A \bar v(\e,\delta,\sigma,\lambda) \right)=0.
\]
We define $\mathbf{A}=\left(\Re(A),\Im(A)\right)$, with $\mathbf{A}\in \R^2$, such that the above equation is reduced into a two dimensional system of the form
\be
\mathbf{G}(\e,\delta,\sigma,\lambda) \mathbf{A}=0.
\label{eqSystem2d}
\ee
By construction and the result from Theorem \ref{thmsolper}, we know that
\[ \mathbf{G}(\e,\delta,0,0) = \left(\begin{matrix} \mathbf{G}_{00}(\e,\delta) &0\\0& 0\end{matrix} \right), \]
with 
\[ \mathbf{G}_{00}(\e,\delta) = 2\left( \hat\phi(k_c)\epsilon^2+\left(1+\frac{\mu_c}{2}\hat\phi''(k_c)\right)\delta^2 \right) + \rmO\left(\left( \hat\phi(k_c)\epsilon^2+\left(1+\frac{\mu_c}{2}\hat\phi''(k_c)\right)\delta^2\right)^{3/2} \right). \]
This implies that periodic solutions $\mathbf{v}_{\e,\delta}$ are neutrally stable with respect to perturbations of the same period $2\pi / k$ as we have a zero eigenvalue associated to the translation invariance of the problem and the other eigenvalue $ \mathbf{G}_{00}(\e,\delta)$ is negative. This gives the proof of the first assertion of Theorem \ref{thmstability}. We now study how the critical $\lambda_c=0$ eigenvalue is perturbed when $0<\sigma\ll 1$.

\begin{lemma}
The critical $\lambda_c=0$ eigenvalue has the following expansion
\[ \lambda_c = \mathbf{g}(\e,\delta) \sigma^2+\rmO(\sigma^4), \text{ as } \sigma \rightarrow 0. \]
To leading order, the coefficient $\mathbf{g}(\e,\delta)$ is given by
\[ \mathbf{g}(\e,\delta) \approx -\frac{k_c^2}{\mathbf{G}_{00}(\e,\delta)}\left( 4\left(1+\frac{\mu_c}{2}\hat\phi''(k_c) \right)\delta^2  +\mathbf{G}_{00}(\e,\delta)\right)\left(1+\frac{\mu_c}{2}\hat\phi''(k_c) \right).\]
% +\rmO\left( \left( \hat\phi(k_c)\epsilon^2+\left(1+\frac{\mu_c}{2}\hat\phi''(k_c)\right)\delta^2\right)^{3/2}\right).\] 
There is an instability whenever $\mathbf{g}(\e,\delta)$ is positive.
\end{lemma}
\begin{Proof}
We want to obtain an expansion of the eigenvalues of the matrix $\mathbf{G}(\e,\delta,\sigma,\lambda)$. We first recall that 
\[
v(\e,\delta,\sigma,\lambda)= D_{A}\Psi(\sqrt{\Gamma},\sqrt{\Gamma},\e,\delta) +i\sigma B(\e,\delta,\sqrt{\Gamma})+\Gamma \rmO(\sigma^2+|\lambda|),
\]
where $\Gamma= \frac{\hat\phi(k_c)\epsilon^2+\left(1+\frac{\mu_c}{2}\hat\phi''(k_c)\right)\delta^2}{\omega}>0$, $D_{A}\Psi(\sqrt{\Gamma},\sqrt{\Gamma},\e,\delta)=\rmO(\sqrt{\Gamma})$ and $|B(\e,\delta,\sqrt{\Gamma})|=\rmO(\sqrt{\Gamma})$. Then, inserting this expansion into \eqref{eqSystem2d}, one obtains 
\[
\mathbf{G}(\e,\delta,\sigma,\lambda)=\left(
\begin{matrix} 
\mathbf{G}_{00}(\e,\delta)+\rho-\lambda & i \beta \\
-i\beta & \rho -\lambda
\end{matrix}
\right)+\Gamma^2\left( \begin{matrix} \rmO(\sigma^2 + |\lambda|) & \rmO(|\sigma|+|\lambda|)  \\  \rmO(|\sigma|+|\lambda|)  & \rmO(\sigma^2 + |\lambda|) \end{matrix} \right),
\]
with $\rho=-k_c^2\left(1+\frac{\mu_c}{2}\hat\phi''(k_c) \right)\sigma^2+\rmO\left(\sigma^2|\delta| \right)$ and $\beta=\left(-2k^2-\mu_c k\hat\phi'(k) \right)\sigma$.

The determinant of $\mathbf{G}(\e,\delta,\sigma,\lambda)$ has an expansion of the form
\[ \det \mathbf{G}(\e,\delta,\sigma,\lambda) = \mu_0+\mu_1\lambda +\rmO(|\lambda|^2), \]
where the coefficients are given by
\begin{align*}
\mu_0(\e,\delta,\sigma)&:= \rho^2 +\mathbf{G}_{00}(\e,\delta)\rho -\beta^2 +\rmO\left( \Gamma^2\sigma^4 \right),\\
\mu_1(\e,\delta,\sigma)&:=-\mathbf{G}_{00}(\e,\delta)-2\rho + \rmO\left( \Gamma^2\sigma^2 \right). 
\end{align*}
We directly note that $\mu_0 \ll \mu_1$, and the critical eigenvalue $\lambda_c$ has thus the expansion $\lambda_c=-\mu_0/ \mu_1+\rmO(\mu_0^2)$. To obtain the expansion of $\mathbf{g}(\e,\delta)$, one further notes that
\[\beta =-2k_c\left(1+\frac{\mu_c}{2}\hat\phi''(k_c) \right)\sigma \delta + \rmO(|\sigma||\delta|^2).\]
\end{Proof}

\begin{Proof}{[of Theorem \ref{thmstability}]}
From the leading order terms in the expansion of $\mathbf{g}(\e,\delta)$, we see that the instability condition $\mathbf{g}(\e,\delta)>0$ is equivalent to
\[ -4\left(1+\frac{\mu_c}{2}\hat\phi''(k_c) \right)\delta^2  -\mathbf{G}_{00}(\e,\delta) <0,\] 
as $\mathbf{G}_{00}(\e,\delta)<0$. The above condition can rewritten as
\[ \frac{-\hat \phi(k_c)}{3\left(1+\frac{\mu_c}{2}\hat\phi''(k_c)\right)} \e^2 < \delta^2< \frac{-\hat \phi(k_c)}{1+\frac{\mu_c}{2}\hat\phi''(k_c)}\e^2,\]
which concludes the second part of the proof of Theorem \ref{thmstability}.
\end{Proof}

\section{Approximate description of modulated fronts using amplitude equation}\label{secAmpEq}

In this section, we will derive an amplitude equation that describes the dynamics of the modulated fronts of (\ref{eq:NLKPP}) for values of $|\mu-\mu_c|\ll 1$.  The amplitude equation is found via a multiple scale analysis of (\ref{eq:NLKPP}) and gives a formal description in the asymptotic limit $\e\to 0$ of the slow modulation of periodic solutions in space and time, see for example \cite{Cross,Mie2}.  This formal calculation will suggest the existence of modulated traveling fronts that we will prove the existence of in section~\ref{secModTF}.

We will use the specific kernel \eqref{eq:kernel} and the associated system \eqref{systemPDE}. We note that the homogeneous stationary state $u=1$ of \eqref{eq:kernel} transforms into $(u,v,w)=(1,3,-2)$ in the system \eqref{systemPDE}, which we transform to the origin to find the system (abusing notation with the same variables),
\begin{subequations}
\begin{align}
u_t &= u_{xx}-\mu (1+u) (v+w)  \\
0 &= v_{xx}-a^2 v+3a^2u  \\
0&= w_{xx}-w-2u. 
\end{align}
\label{eq:systemnearzero}
\end{subequations}
To find an amplitude equation, we write $\e^2=\mu-\mu_c$ and seek solutions of the form,
\be U(t,x)=\e A(X,T)e^{ik_cx}E+c.c.,\label{eq:LOansatz} \ee
where $U=(u,v,w)^T$, $T=\e^2t$ and $X=\e x$.  We re-write (\ref{eq:systemnearzero}) as
\be D U_t=U_{xx}+M_c U+\e^2 M_rU-(\mu_c+\e^2)N(U),\label{eq:condensed}\ee
where $D(1,1)=1$ and all other entries are zero and
\[ M_c=\left(\begin{array}{ccc} 0 & -\mu_c&-\mu_c \\3a^2 & -a^2 & 0 \\ -2 & 0 &-1 \end{array}\right),\quad M_r=\left(\begin{array}{ccc} 0 & -1&-1 \\0 & 0 & 0 \\ 0 & 0 &0 \end{array}\right),\quad N(U)=\left(\begin{array}{c} u(v+w) \\ 0 \\ 0\end{array}\right). \]
We will also need the matrix $L_j=-(jk_c)^2I+M_c$.  The matrix $L_1$ has a kernel which we denote by $E$ that corresponds to a solution $U=Ee^{ik_cx}$ to the differential equation $U_{xx}+M_cU=0$.  We have,
\[ E=(1,\frac{3a^2}{a^2+k_c^2},-\frac{2}{1+k_c^2})^T=(1,\hat{\phi_v}(k_c),\hat{\phi_w}(k_c))^T.\]
We also note
\[ \ \mathrm{coker}(L_1)=\mathrm{span}\left\{ \left(1,-\frac{\mu_c}{a^2+k_c^2},-\frac{\mu_c}{1+k_c^2}\right)^T\right\}.\]
For future reference, we will let $P_1$ denote the above vector in the cokernel of $L_1$.  Note that $d_\nu(0,k_c,\mu_c)=0$ implies that $E$ and $P_1$ are orthogonal. 
  Assuming a leading order ansatz (\ref{eq:LOansatz}),   the higher order and nonlinear terms will generate solutions at order $\e^2$ of the form,
\[ -\mu_cN(Ae^{ik_cx}E+\bar{A}e^{-ik_cx}E).\]  
Thus, we require terms at $\rmO(\e^2)$ that account for these influences.   In particular, we let
\[ E_2=L_2^{-1}N(E), \quad E_0=L_0^{-1}N(E).\]
By our assumption that the mode at $k_c$ is critical at $\e=0$, we have the invertibility of $L_2$ and $L_0$ and we find the formulas,
\begin{eqnarray}
E_2&=& \frac{-\hat{\phi}(k_c)}{4k_c^2+\mu_c\hat{\phi}(2k_c)}\left(1,\hat{\phi}_v(2k_c),\hat{\phi}_w(2k_c)\right)^T\\
E_0&=& -\frac{\hat{\phi}(k_c)}{\mu_c} \left( 1,3, -2\right)^T.
\end{eqnarray}

Thus, our ansatz for the solution of the differential equation is given by,
\be u(t,x)=\e\left(A_1e^{ik_cx}E+c.c.\right)+\e^2\left(A_2e^{2ik_cx}E_2+c.c.\right)+\e^2A_0E_2 +\e^2 \left(A_{1,1}e^{ik_cx}E_1+c.c\right)+\rmO(\e^3),
\ee
where all the amplitudes are functions of the slow space and time variables, i.e. $A_j(X,T)$.  The vector $E_1$ remains to be determined.  

We now plug this into (\ref{eq:condensed}) and solve order by order.
\paragraph{At $\rmO(\e)$}
At leading order we reproduce exactly the linear system,
\[ 0=(A_1e^{ik_cx}+c.c.)L_1E.\]
\paragraph{At $\rmO(\e^2)$}
At next order, we find the equation,
\begin{eqnarray*} 0&=&\left(A_{1,1}e^{ik_cx}L_1E_1+c.c.\right)+\left(2ik_c\frac{\partial A_1}{\partial X}e^{ik_cx}E +c.c.\right) \\
&+& \left(e^{2ik_cx}\left(-\mu_cA_1^2+A_2\right)N(E)+c.c.\right) \\
&+& \left(-2\mu_cA_1\bar{A}_1+A_0\right)N(E) .
\end{eqnarray*}
We eliminate constant terms and those proportional to $e^{\pm 2ik_cx}$ by imposing conditions on $A_{ 2}$ and $A_0$.  We find,
\[ A_2=\mu_cA_1^2, \quad A_0=2\mu_cA_1\bar{A}_1.\]  
We are now left with a linear system of equations describing solutions at $\rmO(\e^2)$ proportional to $e^{ik_cx}$.  Since $E$ and $P_1$ are orthogonal we have that $E\in \mathrm{rng}(L_1)$.  Therefore, we have
\[ A_{1,1}=-2ik_c\frac{\partial A_1}{\partial X}, \quad L_1E_1=E , \quad E_1= \left(1,\frac{3a^2(a^2+k_c^2-1)}{(a^2+k_c^2)^2} , \frac{-2k_c^2}{(1+k_c^2)^2}\right)^T. \]

\paragraph{At $\rmO(\e^3)$}
Continuing to third order in $\e$, we focus only on those terms with a prefactor of $e^{ik_cx}$.  Identifying these terms we find,
\[ \frac{\partial A_1}{\partial T}DE=A_{1,2}L_1E_{1,2}+\frac{\partial^2 A_1}{\partial X^2} E+2ik_c\frac{\partial A_{1,1}}{\partial X}E_1+A_1M_rE-\mu_c^2A_1|A_1|^2 \tilde{E}_2-2\mu_c^2A_1|A_1|^2\tilde{E}_0.\]
We can then write this equation as the vector sum $\gamma P_1+R$, for some $R\in\mathrm{Rg}(L_1)$.  Applying the solvability condition $\gamma=0$ yields the amplitude equation,
\[ \frac{\partial A_1}{\partial T}=\frac{E\cdot P_1}{(DE)\cdot P_1}\frac{\partial^2 A_1}{\partial X^2}+2ik_c\frac{\partial A_{1,1}}{\partial X}\frac{E_1\cdot P_1}{DE\cdot P_1}+A_1\frac{M_RE\cdot P_1}{DE\cdot P_1}-\mu_c^2\frac{(\tilde{E}_2+2\tilde{E}_0)\cdot P_1}{(DE)\cdot P_1}A_1|A_1|^2.\]
Exploiting the formula for $A_{1,1}$, the fact that $E\cdot P_1=0$ and using that $DE\cdot P_1=1$ this immediately simplifies to,
\[ \frac{\partial A_1}{\partial T}=4k_c^2E_1\cdot P_1\frac{\partial^2 A_1}{\partial X^2}+A_1M_RE\cdot P_1-\mu_c^2(\tilde{E}_2+2\tilde{E}_0)\cdot P_1A_1|A_1|^2.\]
Finally, we have,
\begin{eqnarray*}
E_1\cdot{P_1}&=&  1-\mu_c \left( \frac{3a^2(a^2+k_c^2-1)}{(a^2+k_c^2)^3} - \frac{2k_c^2}{(1+k_c^2)^3} \right) \\
M_rE\cdot P_1&=& -\hat{\phi}(k_c)\\
\tilde{E}_2&=& N(E_2+E)-N(E_2)-N(E) \\
&=& \frac{-\hat{\phi}(k_c)}{4k_c^2+\mu_c\hat{\phi}(2k_c)}\left(\hat{\phi}(k_c)+\hat{\phi}(2k_c),0,0\right)^T \\
\tilde{E}_0&=& N(E_0+E)-N(E_0)-N(E) \\
&=& \left(-\frac{\hat{\phi}(k_c)+\hat{\phi}(k_c)^2}{\mu_c},0,0\right)^T.
\end{eqnarray*}
Putting this all together, we find the amplitude equation,
\be \frac{\partial A_1}{\partial T}=\zeta\frac{\partial^2 A_1}{\partial X^2}-\hat{\phi}(k_c)A_1+\omega A_1|A_1|^2,\label{eqAmpMod}\ee
where 
\be \zeta := 4k_c^2 \mu_c \left( \frac{3a^2}{(a^2+k_c^2)^3} - \frac{2}{(1+k_c^2)^3} \right)=1+\frac{\mu_c}{2}\hat{\phi}''(k_c), \label{eqcoeffzeta}\ee
and $\omega$ is given in equation \eqref{eqcoeffomega}. Traveling wave solutions of \eqref{eqAmpMod}, and in particular traveling fronts $A_1(T,X)=A_1(X-sT)=\widetilde{A}(Z)$ satisfy
\be \zeta\frac{\partial^2 \widetilde{A}}{\partial Z^2}+s\frac{\partial \widetilde{A}}{\partial Z}-\hat{\phi}(k_c)\widetilde{A}+\omega \widetilde{A}|\widetilde{A}|^2=0,\label{eqAmpModTF}
\ee
where $Z=X-sT$. This equation is invariant under the phase-shift transformation $\widetilde{A}\rightarrow \widetilde{A}e^{i\theta}$, $\theta\in\R$. Then restricting ourselves to real $A_1\in\R$, one can easily prove, using phase-plane analysis, the existence of a heteroclinic connection between $(0,0)$ and $\left(0, \sqrt{\frac{\hat\phi(k_c)}{\omega}}\right)$. And, when $s^2\geq-4\hat\phi(k_c)\zeta$ this front is monotone. The amplitude equation leads one to expect that 
\[ u(t,x) = \e \widetilde{A}(X-sT) \cos(k_cx), \]
where $\widetilde{A}$ is solution of  \eqref{eqAmpModTF}, is a good approximation of the modulated traveling fronts of \eqref{eq:NLKPP}.  We prove the existence of such solutions in the following section.

\section{Modulated Traveling Fronts}\label{secModTF}

In this section, we will construct modulated traveling front solutions for $\mu\approx \mu_c$ via  center manifold reduction.  We first set up the problem.  Then in section~\ref{sec:linear} we study the spectrum of the linear operator first for $\e=0$ and then for $\e>0$.  In section~\ref{sec:CMR} we apply a version of the Center Manifold Theorem due to Haragus and Schneider \cite{Har-Sch} and compute the reduced equations on the center manifold.  After several changes of coordinates, we observe in \ref{sec:IDamplitude} that the leading order dynamics on the center manifold are equivalent to the formal leading order amplitude equations.  Finally, in section~\ref{sec:heteroclinic} we study the equations on the center manifold and find heteroclinic connections that correspond to the modulated traveling front solution.

We seek solutions of (\ref{eq:condensed}) of the form,
\[ U(t,x)=W(x-\e st,x)=\sum_{n\in\mathbb{Z}} W_n(x-\e st)e^{-ink_cx}.\]
Let $\xi=x-\e st$.  Then, plugging this ansatz into (\ref{eq:condensed}), we find a system of six coupled ordinary differential equations for each Fourier mode $n$.  Let 
\[ X := \left( X_n\right)_{n\in \Z}, \quad X_n :=\left(W_n^u,\partial_\xi W_n^u,W_n^v,\partial_\xi W_n^v,W_n^w,\partial_\xi W_n^w\right)^T.\]
Then, (\ref{eq:condensed}) is transformed to a system of equations,
\be \partial_\xi X_n=\mathcal{M}_n^\epsilon X_n-(\mu_c+\e^2)\mathcal{N}_n(X), \quad n\in \Z. \label{eqTWmode}\ee
The matrix $\mathcal{M}_n^\e$ is defined as
\be
\mathcal{M}_n^\epsilon:=\left(\begin{matrix} 0 & 1 & 0&0 & 0 & 0 \\
n^2k_c^2 & 2 i n k_c-\e s & \mu_c+\e^2 & 0 & \mu_c+\e^2 & 0 \\
0 & 0 & 0&1 & 0 & 0 \\
-3a^2 & 0 & a^2+n^2k_c^2 & 2 i n k_c & 0 & 0 \\
0 & 0 & 0&0 & 0 & 1 \\
2 & 0 & 0 & 0 & 1+n^2k_c^2 & 2 i n k_c \end{matrix}\right),\label{eqMatrixMn}\ee
and the nonlinear term as
\be
\mathcal{N}_n(X):=\left(\begin{matrix}
0 \\
\underset{p+q=n}{\sum}\left(W_p^uW_q^v+W_p^uW_q^w \right) \\
0 \\ 0 \\ 0 \\ 0
\end{matrix} \right).\label{eqNonlinear}
\ee

We recall that neutrality of the critical mode gives the identity (see (\ref{eq:disp})),
\be -k_c^2-\mu_c\left( \frac{3a^2}{a^2+k_c^2}-\frac{2}{1+k_c^2}\right)=0.\label{eq:crit1}\ee
Furthermore, the criticality of the neutral mode implies the identity,
\be 1-\mu_c\left( \frac{3a^2}{(a^2+k_c^2)^2}-\frac{2}{(1+k_c^2)^2}\right)=0.\label{eq:crit2}\ee

\subsection{The linear operator $\mathcal{M}_n^\e$}\label{sec:linear}

In order to capture the nature of the spectrum of $\mathcal{M}_n^\epsilon$, for $n\in \mathbb{Z}$ and $0<\epsilon \ll 1$, we first set $\e=0$. One can easily check that the characteristic polynomial of $\mathcal{M}_n^0$ simplifies into
\be p_n(\lambda) = (\lambda^2-2ink_c\lambda -1-a^2-n^2k_c^2-2k_c^2 )\left(\lambda-i(n-1)k_c \right)^2\left(\lambda-i(n+1)k_c \right)^2, \quad n\in \mathbb{Z}.\label{eqCPn} \ee
This implies that the spectrum $\sigma_n^0$ of $\mathcal{M}_n^0$ is
\[ \sigma_n^0=\left\{ i(n\pm1)k_c, ink_c\pm \sqrt{1+a^2+2k_c^2} \right\}, \quad n\in \mathbb{Z}. \]
Each eigenvalue $\lambda_{\pm,n}=i(n\pm1)k_c$ has algebraic multiplicity two and geometric multiplicity one. We also define $\beta_{\pm,n}=ink_c\pm \sqrt{1+a^2+2k_c^2}$. 
Of interest for our forthcoming computations are the eigenvector and generalized eigenvector associated to $\lambda_{+,-1}$ and $\lambda_{-,1}$:
\[ \left(1,0,\hat\phi_v(k_c),0,\hat\phi_w(k_c),0\right)^T, \quad \left( 0, 1, -\frac{2ik_c}{a^2+k_c^2}\hat\phi_v(k_c),\hat \phi_v(k_c),-\frac{2ik_c}{1+k_c^2}\hat\phi_w(k_c),\hat\phi_w(k_c)\right)^T. \]
We now study how these eigenvalues are perturbed away when we turn on the parameter $\epsilon$. We discuss several cases.

\begin{figure}%[h!]
\centering
\includegraphics[width=0.5\linewidth]{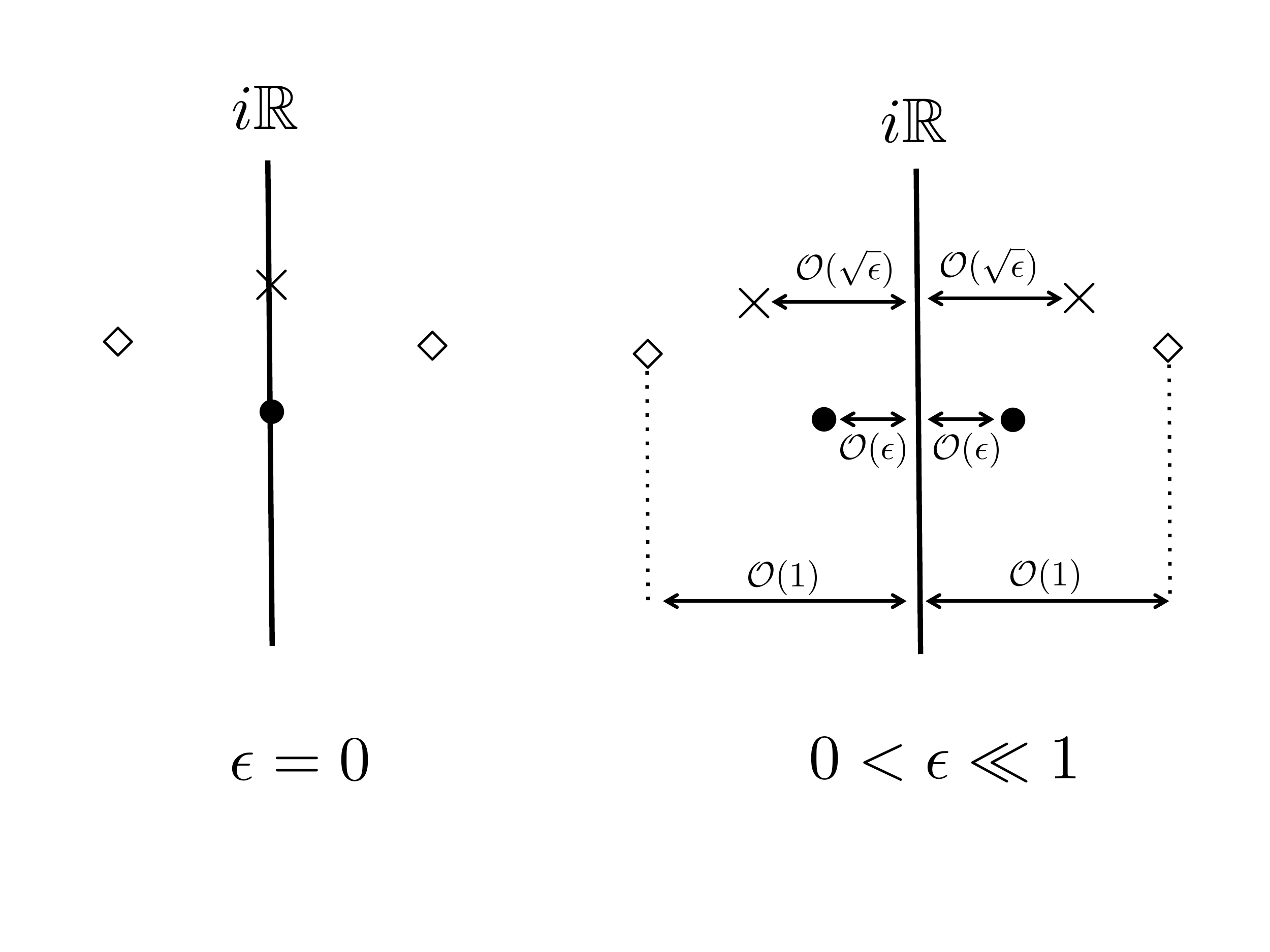}
\caption{Illustration of the splitting of the eigenvalues of $\mathcal{M}_1^\e$ for $0<\e<\ll1$. The dot represents the double $0$ eigenvalue, the cross represents the double $2ik_c$ eigenvalue. The diamond represents the two eigenvalue $ik_c\pm\sqrt{1+a^2+2k_c^2}$.}
\label{fig:splitting}
\end{figure}

\paragraph{Case $\lambda_{\pm,\mp 1}=0$} We only deal with $\lambda_{-,1}$ as the computations for the other eigenvalue follow along the same lines. We first remark that the characteristic polynomial associated to $\mathcal{M}_1^\epsilon$ can be written as follows
\[ p^\epsilon_1(\lambda) = p_1(\lambda) + \epsilon \mathcal{R}_1^\epsilon(\lambda), \]
where $p_1(\lambda)$ is given in equation \eqref{eqCPn} and $\mathcal{R}_1^\epsilon(\lambda)$ is of order five in $\lambda$, which admits the Taylor expansion
\[ \mathcal{R}_1^\epsilon(\lambda)  = \epsilon \alpha_0 + \alpha_1 \lambda + \mathcal{O}\left( \epsilon|\lambda|+|\lambda|^2 \right),\]
as $\epsilon \rightarrow 0$ and $|\lambda| \rightarrow 0$. Here, we have 
\be \alpha_0 := 2k_c^2-3a^2k_c^2-a^2 \text{ and } \alpha_1 := (1+k_c^2)(a^2+k_c^2)s.\label{eqAlpha01}\ee
As a consequence, we can look for eigenvalues $\lambda_{-,1}^\epsilon$ which can be expanded as
\[ \lambda_{-,1}^\epsilon = \epsilon \chi + \mathcal{O}\left(\epsilon^2\right), \]
where $\chi\in \mathbb{C}$ can be computed explicitly. Looking at terms of order $\mathcal{O}(\epsilon^2)$ in $p^\epsilon_1(\lambda_{-,1}^\epsilon)=0$, we find that $\chi$ satisfies the quadratic equation
\[ \alpha_2\chi^2+\alpha_1\chi+\alpha_0=0, \quad \alpha_2:=4k_c^2(1+a^2+3k_c^2).\]
This equation has roots,
\[ \chi_\pm = \frac{-\alpha_1\pm \sqrt{\Delta}}{2\alpha_2}, \]
where $\Delta:=\alpha_1^2-4\alpha_0\alpha_2$. Here, we impose $\Delta>0$ which gives a condition on $s$. As a conclusion, the algebraically double eigenvalue $\lambda_{-,1}$ perturbs into two eigenvalues with asymptotics:
\[ \lambda_{-,1}^{\epsilon,\pm} = \epsilon \chi_\pm + \mathcal{O}\left(\epsilon^2\right), \text{ as } \epsilon \rightarrow 0. \]
One can then easily check that $\lambda_{+,-1}^{\epsilon,\pm}$ satisfies the same leading order expansion. 

\paragraph{Case $\lambda_{\pm,n}$, $n\neq \pm 1$} Reproducing the same analysis, we look for eigenvalues $\lambda_{\pm,n}^\epsilon$ solutions of
\[ p^\epsilon_n(\lambda) = p_n(\lambda) + \epsilon \mathcal{R}_n^\epsilon(\lambda),\]
where $p_n(\lambda)$ is given in equation \eqref{eqCPn} and $\mathcal{R}_n^\epsilon(\lambda)$ is of order five in $\lambda$. Looking for eigenvalue $\lambda_{\pm,n}^\epsilon$ which can be expanded as
\[ \lambda_{\pm,n}^\epsilon = i(n\pm1)k_c+\sqrt{\epsilon} z + \mathcal{O}\left(\epsilon\right), \]
where $z\in \mathbb{C}$, we obtain, at order $\mathcal{O}(\epsilon)$ in  $p^\epsilon_n(\lambda_{\pm,n}^\epsilon)=0$, the compatibility condition
\[\alpha_2z^2+i\alpha_1(n\pm1)k_c=0. \]
This gives complex solutions of the form
\[ z_\pm = \pm \sqrt{\frac{-i\alpha_1(n\pm1)k_c}{\alpha_2}}, \quad n \in \mathbb{Z}, \quad n\neq \pm1.\]

\paragraph{Case $\beta_{\pm,n}$} For any $n\in \mathbb{Z}$, $\Re(\beta_{\pm,n})=\pm \sqrt{1+a^2+2k_c^2}\neq 0$ and as $0<\epsilon \ll 1$, these eigenvalues will be perturbed into $\beta_{\pm,n}^\epsilon$ with $\Re(\beta_{\pm,n}^\epsilon)=\rmO(1)$ by hyperbolicity as $\e\rightarrow 0$.

\paragraph{Conclusion} This perturbation analysis shows that eigenvalues $\lambda_{\pm,\mp1}^{\epsilon,\pm}$ are "more central" than all other eigenvalues, as $\Re(\lambda_{\pm,\mp1}^{\epsilon,\pm}) = \mathcal{O}(\epsilon)$, see Figure \ref{fig:splitting}. These eigenvalues are an $\mathcal{O}(1)$ distance from all the other eigenvalues and therefore the spectral projection onto the center spectrum is uniformly bounded in $\e$.  More precisely, we have obtained the following proposition.

\begin{proposition}
For each fixed $s\in\R$ such that $\Delta>0$, there are positive constants $d_0$, $d_1$ and $\e_0$ such that the following is true for all $\e\in (0,\e_0)$: there are precisely four eigenvalues within $\underset{n\in\Z}{\cup}\sigma_n^\e$ counted with multiplicity within the strip $|\Re(\lambda)|\leq d_0 \e$, while the remainder of the eigenvalues $\lambda$ satisfies $|\Re(\lambda)|\geq d_1 \sqrt{\e}$. Furthermore, the spectral projection $\mathcal{P}_c$ onto the $\e$-dependent, four-dimensional eigenspace associated with the eigenvalues within distance $\rmO(\e)$ of the imaginary axis is bounded uniformly in $\e$.
\end{proposition}

\subsection{Center manifold reduction}\label{sec:CMR}

In this section, we revisit a center manifold result developed by Eckmann and Wayne in \cite{Eck-Way} and adapt it to our problem along the lines of Haragus and Schneider \cite{Har-Sch}.

\subsubsection{Functional setup}

The function spaces that we will use are analogous to those used in \cite{Eck-Way,Har-Sch} and we borrow their notation. We denote $\E_0$ the direct sum $\overset{\infty}{\underset{n=0}{\bigoplus}} \C^6$, and if $X\in \E_0$, we denote by $X_{n,j}$, $n\geq0$ and $j=1,\dots,6$, the components of $X$ in the decomposition of $\E_0$. Furthermore, let $\E$ be the subset of $\E_0$ in which $X_{0,j}$, $j=1,\dots,6$ are real. We have a one-to-one map $I_\xi$, $\xi\in\R$ from the space of real, sufficiently differentiable functions of form
\[ W(\xi,x)=\sum_{n\in\Z}W_n(\xi)e^{-ink_cx}, \quad W_n(\xi)=\left(W_n^u(\xi),W_n^v(\xi),W_n^w(\xi) \right)^T, \]
into $\E$. It is defined by
\be I_\xi(W)_n=\left(W_n^u,\partial_z W_n^u|_{z=\xi},W_n^v,\partial_z W_n^v|_{z=\xi},W_n^w,\partial_z W_n^w|_{z=\xi}\right)^T \longleftrightarrow X_n. \label{definitionIxi}\ee
We remark that $W_{-n}=\bar W_n$, so that $I_\xi(W)_n\in\C^6$, $n>0$ and $I_\xi(W)_0\in\R$ uniquely determine $W$. We define an inner product on $\E_0$ 
\[ \langle X,Y \rangle_\ell = \sum_{n = 0}^\infty (1+n^2)^\ell \langle X_n , Y_n \rangle_{\C^6} \]
and define the Hilbert space $\mathcal{H}^\ell_\C(\E_0)=\left\{ X\in \E_0 ~|~ \langle X,X \rangle_\ell <\infty \right\}.$ The nonlinearity $\mathcal{N}:\E_0\rightarrow \E_0$ defined as
\[ \mathcal{N}_n(X)=\left( 0,\underset{p+q=n}{\sum}\left(X_{p,1}X_{q,3}+X_{p,1}X_{q,5} \right), 0,0,0,0\right)^T, \quad n\geq0, \]
is a continuously differentiable function from $\mathcal{H}^\ell_\C(\E_0)$ into $\mathcal{H}^\ell_\C(\E_0)$ as long as $\ell>1/2$. Using the fact that $\mathcal{H}^\ell_\C(\E_0)$ is a Banach algebra, we further have the estimate:
\[ \| \mathcal{N}(X) \|_\ell \leq C \| X \|_l^2, \quad X\in \mathcal{H}^\ell_\C(\E_0), \]
for some constant $C>0$ depending on $\ell$ and the nonlinearity. Finally, we define the bilinear map $\widetilde{\mathcal{N}}: \E_0\times\E_0\rightarrow \E_0$ with components 
\be \widetilde{\mathcal{N}}_n(X,Y)=\frac{1}{2}\left( 0,\underset{p+q=n}{\sum}\left(X_{p,1}Y_{q,3}+Y_{p,1}X_{q,3}+X_{p,1}Y_{q,5}+Y_{p,1}X_{q,5}  \right), 0,0,0,0\right)^T, \quad n\geq0,\label{eqbilinear}\ee
for any $(X,Y)\in \E_0\times\E_0$. Note that in the definition of $\mathcal{N}_n(X)$, the sum over the couples $p+q=n$, $(p,q)\in \Z^2$ has to be understood as follows. If, for example  $p<0$, then the corresponding term $X_{p,j}$ is replaced by $\bar X_{-p,j}$.

\subsubsection{Reduction to a center manifold}

We can now apply a center manifold reduction to our infinite-dimensional system \eqref{eqTWmode}, which takes the form
\be\partial_\xi X=\mathcal{M}^\epsilon X-(\mu_c+\e^2)\mathcal{N}(X), \quad X\in \E,\label{eqformalTW}\ee
where, component wise, we have $\left(\mathcal{M}^\epsilon X\right)_n=\mathcal{M}^\epsilon_n X_n$ and $\left(\mathcal{N}(X)\right)_n=\mathcal{N}_n(X)$ for $n\in\N$ (see \eqref{eqMatrixMn} and \eqref{eqNonlinear}). Here, the center directions will consist of the subspace of $\E$ corresponding to the eigenvalues whose real part are of order $\rmO(\e)$. From the analysis of the previous section, we deduce that there are only two center directions corresponding to the eigenvalues $\lambda_\pm^\e:=\lambda_{-,1}^{\e,\pm}=\e \chi_\pm+\rmO(\e^2)$ of $\mathcal{M}_1^\e$.  We denote by $\varphi_\pm^\e\in \C^6$ the two corresponding eigenvectors that satisfy
\[ \mathcal{M}_1^\e \varphi_\pm^\e = \lambda_\pm^\e \varphi_\pm^\e. \]
Furthermore, we denote by $\psi_\pm^\e\in \C^6$ the eigenvectors of the adjoint matrix $\left( \mathcal{M}_1^\e\right)^*$ of $\mathcal{M}_1^\e$ corresponding to the eigenvalue $\bar{\lambda}_\pm^\e$. An expansion of these eigenvectors can be found in the Appendix \ref{AppIP} and some simple computations show that
\[
\langle \psi_\pm^\e, \varphi_\pm^\e \rangle =  \pm \frac{\sqrt{\Delta}}{\alpha_2}\e \text{ and } \langle \psi_\mp^\e, \varphi_\pm^\e \rangle =0.
\]
We also denote by $\Phi_\pm^\e \in \E$ (respectively $\Psi_\pm^\e\in\E$) the vector in $\E$ defined as
\[ \left( \Phi_\pm^\e \right)_1=\varphi_\pm^\e, \quad  \left( \Phi_\pm^\e \right)_n= 0_{\C^6}, \quad n\neq 1.\]
We can now define the spectral projection $\mathcal{P}_c^\e:\E\rightarrow \E_c$ where $\E_c=\left\{\Phi_+^\e,\Phi_-^\e \right\}$ via
\[ \mathcal{P}_c^\e X= c_+^\e \langle \Psi_+^\e,X\rangle_\ell \Phi_+^\e+c_-^\e \langle \Psi_-^\e,X\rangle_\ell \Phi_-^\e,  \]
where $c_\pm^\e$ are normalization constant such that $c_\pm \langle \Psi_\pm^\e,\Phi_\pm\rangle_\ell=1$. Noting that $\langle \Psi_\pm^\e,\Phi_\pm\rangle_\ell = 2^\ell \langle \psi_\pm^\e, \varphi_\pm^\e \rangle$, we have the asymptotics
\[ c_\pm^\e = \frac{\pm \alpha_2}{2^\ell \e\sqrt{\Delta}}+\rmO(1), \text{ as } \e\rightarrow 0. \]
Finally, we have shown in the previous section that there exists a constant $d_1>0$ such that the spectrum of $\mathcal{M}^\e|_{(\text{id}-\mathcal{P}_c^\e)\E}$ satisfies the estimate
\[ \sigma\left( \mathcal{M}^\e|_{(\text{id}-\mathcal{P}_c^\e)\E} \right)\subset \left\{ \lambda \in \C ~|~ |\Re(\lambda)| \geq d_1 \sqrt{\e} \right\}. \]

We can now apply the center manifold result of Haragus and Schneider \cite{Har-Sch} (see also \cite{Har-Ioo}, \S2 Remark 3.6) which is a generalization to quadratic nonlinearities of the center manifold result initially developed by Eckmann and Wayne \cite{Eck-Way} for the cubic case. The crucial step of the proof of \cite{Har-Sch} is to show that the product of the Lipschitz constant of the nonlinearity, restricted to a ball of radius $\e^{2/3+\gamma}$, with $0<\gamma<1/3$, centered at the rest state $X=0$, and the inverse of the spectral gap between the center and the hyperbolic part of $\mathcal{M}^\e$ is smaller than one for $\e>0$. To obtain such a control on the nonlinear terms, the key idea is to decompose the solutions of \eqref{eqformalTW} into a central part along the projection $\mathcal{P}_c^\e$ and a hyperbolic part which is along the rest of the spectrum. This hyperbolic part is then further reduced into a part along the eigenvalues of order $\rmO\left(\e^{1/2}\right)$ and a part along eigenvalues of order $\rmO(1)$. With this decomposition in hand, Haragus and Schneider \cite{Har-Sch} have shown that provided $\mathcal{P}_c^\e \widetilde{\mathcal{N}}(\mathcal{P}_c^\e X,\mathcal{P}_c^\e X)=0$, which is our case here, a normal form transformation, which eliminates the potentially dangerous quadratic terms, leads to the following center manifold theorem; see Theorem 6.3 in \cite{Har-Sch}.

\begin{proposition}\label{prop:cm}
For $0<\gamma<1/3$ and $\e>0$ sufficiently small, there exists neighborhoods of the origin $\mathcal{U}^\e\subset \E_c$, $\mathcal{V}^\e\subset (\text{id}-\mathcal{P}_c^\e)\E$, and for any $m<\infty$, a $\mathcal{C}^m$-map $\Theta^\e:\mathcal{U}^\e\rightarrow \mathcal{V}^\e$ having the following properties.
\begin{enumerate}
\item[(i)] All bounded solutions of \eqref{eqformalTW} within $\mathcal{U}^\e\times \mathcal{V}^\e$ are on the center manifold, i.e., 
\[ X = X_c+\Theta^\e(X_c). \]
\item[(ii)] The center manifold is tangent to the center eigenspace, i.e.,
\[ \left\| \Theta^\e(X_c) \right\|_\ell = \rmO_\e\left(\left\| X_c \right\|_\ell^2 \right).  \] 
\item[(iii)] The neighborhood $\mathcal{U}^\e$ is of size $\rmO(\e^{2/3+\gamma})$.
\end{enumerate}
\end{proposition}

\subsubsection{Reduced system}

We now compute the reduced equation on the center manifold that is obtained by projecting equation \eqref{eqformalTW} with $\mathcal{P}_c^\e$, and we obtain
\be \frac{d X_c}{d\xi} = \mathcal{M}^\e X_c -(\mu_c+\e^2)\mathcal{P}_c^\e\left(\mathcal{N}(X_c+\Theta^\e(X_x))\right).\label{eqReduced1}\ee
We introduce the coordinate $X_c=x_+ \Phi_+^\e+x_-\Phi_-^\e$ on the center manifold such that
\[ X = x_+ \Phi_+^\e+x_-\Phi_-^\e+\Theta^\e(x_+,x_-), \]
where  $\Theta^\e(x_+,x_-)$ has the expansion
\[\Theta^\e(x_+,x_-)= \sum_{|{\bf m}|=2}x_+^{m_1}x_-^{m_2}\bar x_+^{m_3}\bar x_-^{m_4}\Phi_{{\bf m}}^\e+\rmO_\e\left( |x_++x_-|^3 \right), \]
and $\Phi_{{\bf m}}^\e$ all belong to $(\text{id}-\mathcal{P}_c^\e)\E$. Here ${\bf m}\in \N^4$ is a multi-index.  As a consequence, we obtain a set of two ordinary differential equations given by
\begin{subequations}
\begin{align}
\frac{dx_+}{d\xi}&=\lambda_+^\e x_+-(\mu_c+\e^2)c_+^\e \langle \Psi_+^\e, \mathcal{N}(x_+ \Phi_+^\e+x_-\Phi_-^\e+\Theta^\e(x_+,x_-)) \rangle_\ell , \\
\frac{dx_+}{d\xi}&=\lambda_-^\e x_--\mu_c (\mu_c+\e^2)c_-^\e \langle \Psi_-^\e, \mathcal{N}(x_+ \Phi_+^\e+x_-\Phi_-^\e+\Theta^\e(x_+,x_-)) \rangle_\ell .
\end{align}
\label{eqReduced2}
\end{subequations}
Note also that $\Theta^\e$ satisfies
\be
{\bf D}_X\Theta^\e(X_c) \frac{dX_c}{d\xi} = \mathcal{M}^\e X_c-(\mu_c+\e^2){\mathcal{P}_c^\e}^\perp (\mathcal{N}(X_c+\Theta^\e(X_c))),\label{eqTheta}
\ee
where ${\mathcal{P}_c^\e}^\perp:=\text{id}-\mathcal{P}_c^\e$.
From the definition of $\Psi_\pm^\e$ and $ \mathcal{N}$ we have that
\[  \langle \Psi_\pm^\e, \mathcal{N}(x_+ \Phi_+^\e+x_-\Phi_-^\e+\Theta^\e(x_+,x_-)) \rangle_\ell = 2^\ell \langle \psi_\pm^\e , \mathcal{N}_1(x_+ \Phi_+^\e+x_-\Phi_-^\e+\Theta^\e(x_+,x_-))\rangle. \]
Furthermore, using the identity $\mathcal{N}_1(X)=\widetilde{\mathcal{N}_1}(X,X)$ and the expansion of $\Theta^\e$, we can deduce that 
\[ \mathcal{N}_1(x_+ \Phi_+^\e+x_-\Phi_-^\e+\Theta^\e(x_+,x_-)) = \rmO_\e( |x_++x_-|^3).  \]
Indeed, from the definition of $\widetilde{\mathcal{N}}_1$ in \eqref{eqbilinear} and the fact that $\left( \Phi_\pm^\e \right)_1=\varphi_\pm^\e$ and $\left( \Phi_\pm^\e \right)_n= 0_{\C^6}$ for $n\neq 1$, we have that all the quadratic terms are of the form $\widetilde{\mathcal{N}}_1(\Phi_\pm^\e,\Phi_\pm^\e)=0$. To give an example of the cubic terms appearing in $\mathcal{N}_1(x_+ \Phi_+^\e+x_-\Phi_-^\e+\Theta^\e(x_+,x_-))$, we have that 
\[\widetilde{\mathcal{N}}_1(x_+ \Phi_+^\e,x_+^2\Phi_{2,0}^\e)=\left( 0, \left(\varphi_+^\e \right)^u \left(\left(\Phi_{2,0}^\e \right)^v_0 + \left(\Phi_{2,0}^\e \right)^w_0 \right)x_+^3+\left(\bar{\varphi}_+^\e \right)^u \left(\left(\Phi_{2,0}^\e \right)^v_2 + \left(\Phi_{2,0}^\e \right)^w_2 \right)x_+|x_+|^2 ,0,0,0,0\right)^T.  \]

Since the linearization leaves each sub-system invariant, we may work on each mode. We denote $\Theta^\e_n$ the $n$th mode of $\Theta^\e$. If $x_c:=x_+\varphi_+^\e+x_-\varphi_-^\e$, equation \eqref{eqTheta} can be seen as
\be {\bf D}_x\Theta^\e_n(X_c) \frac{dx_c}{d\xi} = \mathcal{M}^\e_n x_c-(\mu_c+\e^2){\mathcal{P}_c^\e}^\perp (\mathcal{N}_n(X_c+\Theta^\e(X_c))).\label{eqThetaN} \ee

\paragraph{The case $n=0$} We suppose that 
\[ \Theta^\e_0(x_+,x_-)=\sum_{|{\bf m}|=2}x_+^{m_1}x_-^{m_2}\bar x_+^{m_3}\bar x_-^{m_4}\theta_{{\bf m}}^0(\e)+\rmO_\e\left( |x_++x_-|^3 \right), \quad \theta_{{\bf m}}^0(\e)\in \C^6, \]
and we insert this expression into \eqref{eqThetaN} to find a hierarchy of equations in $x_+^{m_1}x_-^{m_2}\bar x_+^{m_3}\bar x_-^{m_4}$ which yields the equation
\[ \Lambda_{{\bf m}} \theta_{{\bf m}}^0(\e) = \mathcal{M}_0^\e \theta_{{\bf m}}^0(\e) -(\mu_c+\e^2)\mathcal{N}_{0,{\bf m}}^\e, \]
where 
\[ \Lambda_{{\bf m}} =\sum_{j=1}^4 \lambda_j m_j, \quad \lambda_{1,2}=\lambda_\pm^\e,\quad  \lambda_{3,4}=\bar\lambda_\pm^\e,\]
and $\mathcal{N}_{0,{\bf m}}$ is the nonlinear remainder part that is of the form
\[ \mathcal{N}_{0,{\bf m}}^\e = (0, {\bf N}_{0,{\bf m}}^\e, 0,0,0,0)^T. \]
More precisely, we have that ${\bf N}_{0,{\bf m}}^\e=0$ for all $|{\bf m}|=2$ except for $(1,0,1,0)$, $(1,0,0,1)$, $(0,1,1,0)$ and $(0,1,0,1)$. For example, we have
\[ {\bf N}_{0,(1,0,1,0)}^\e= \varphi_+^u\left(\bar \varphi_+^v +\bar \varphi_+^w \right)+ \bar\varphi_+^u\left( \varphi_+^v+\varphi_+^w \right)\longrightarrow 2 \hat \phi(k_c),\]
as $\e\rightarrow0$. In fact, one has that ${\bf N}_{0,{\bf m}}^\e\longrightarrow 2 \hat \phi(k_c)$ as $\e\rightarrow0$ in the three other cases. Then we can compute,
\begin{align*}  \theta_{{\bf m}}^0(\e) &= (\mu_c+\e^2)\left(\mathcal{M}_0^\e-\Lambda_{{\bf m}}\text{id}\right)^{-1}\mathcal{N}_{0,{\bf m}}^\e\\
&=\mu_c\left( \mathcal{M}_0^0\right)^{-1} \mathcal{N}_{0,{\bf m}}^0 + \rmO(\e)\\
&={\bf N}_{0,{\bf m}}^0(1,0,3,0,-2,0)^T+ \rmO(\e)\\
&={\bf N}_{0,{\bf m}}^0(1,0,\hat\phi_v(0),0,\hat\phi_w(0),0)^T+ \rmO(\e).
\end{align*}

\paragraph{The case $n=1$} The case $n=1$ is slightly different from the previous one. Indeed, as we already noticed $\mathcal{N}_1(x_+\Phi_+^\e+x_-\Phi_-^\e)=0$, we have equations of the form
\[\Lambda_{{\bf m}} \theta_{{\bf m}}^1(\e) = \mathcal{M}_1^\e \theta_{{\bf m}}^1(\e).\]
We note that $\mathcal{M}_1^0$ is not invertible, but the above equation admits solution of the form
\[  \theta_{{\bf m}}^1(\e) = \left(1,0,\hat\phi_v(k_c),0,\hat\phi_w(k_c),0\right)^T+\rmO(\e). \]

\paragraph{The case $n=2$} The same general picture holds when $n = 2$ as in the case $n=0$. The invertibility of $\mathcal{M}_2^\e$ implies that the quadratic coefficients can be computed using the formula
\[ \Lambda_{{\bf m}} \theta_{{\bf m}}^2(\e) = \mathcal{M}_2^\e \theta_{{\bf m}}^2(\e) -(\mu_c+\e^2)\mathcal{N}_{2,{\bf m}}^\e. \]
The main difference here is the vector $\mathcal{N}_{2,{\bf m}}^\e$ which is now nonzero only for ${\bf m}$ equal to $(2,0,0,0)$, $(1,1,0,0)$ or $(0, 2, 0, 0)$. We find the leading order expansion
\begin{align*}  \theta_{{\bf m}}^2(\e) &= (\mu_c+\e^2)\left(\mathcal{M}_2^\e-\Lambda_{{\bf m}}\text{id}\right)^{-1}\mathcal{N}_{2,{\bf m}}^\e\\
&=\mu_c\left( \mathcal{M}_2^0\right)^{-1} \mathcal{N}_{2,{\bf m}}^0 + \rmO(\e)\\
&=\frac{\mu_c}{4k_c^2+\mu_c\hat \phi(2k_c)}{\bf N}_{2,{\bf m}}^0\left(1,0,\frac{3a^2}{a^2+4k_c^2},0,\frac{-2}{1+4k_c^2},0\right)^T+ \rmO(\e).
\end{align*}
One can check that ${\bf N}_{2,{\bf m}}^0=\hat\phi(k_c)$.

\paragraph{Conclusion} From the above case study, one can easily check that in the inner product $\langle \psi_\pm^\e , \mathcal{N}_1(x_+ \Phi_+^\e+x_-\Phi_-^\e+\Theta^\e(x_+,x_-))\rangle$, as $\e\rightarrow0$. At cubic order, one finds
\[ \langle \psi_\pm^\e , \mathcal{N}_1(x_+ \Phi_+^\e+x_-\Phi_-^\e+\Theta^\e(x_+,x_-))\rangle =-\frac{\omega(1+k_c^2)}{\mu_c^2\kappa_0}(x_++x_-)|x_++x_-|^2,\]
% -\frac{\hat\phi(k_c)(1+k_c^2)}{\mu_c\kappa_0}\left( \frac{\mu_c(\hat{\phi}(k_c)+\hat{\phi}(2k_c))}{4k_c^2+\mu_c\hat{\phi}(2k_c)}+2(1+\hat{\phi}(k_c)) \right)(x_++x_-)|x_++x_-|^2 \\
%&= -\frac{\omega(1+k_c^2)}{\mu_c^2\kappa_0}(x_++x_-)|x_++x_-|^2, \end{align*}
where $\omega$ is defined in equation \eqref{eqcoeffomega} and $\kappa_0$ is defined in (\ref{eq:kappa0}). Putting everything together, we obtain reduced equations for the flow on the center manifold of the form

\begin{subequations}
\begin{align}
\frac{dx_+}{d\xi}&=\lambda_+^\e x_++\frac{\omega c_+^\e 2^\ell(1+k_c^2)}{\mu_c\kappa_0}(x_++x_-)|x_++x_-|^2 + \rmO_\e\left( |x_++x_-|^4\right), \\
\frac{dx_+}{d\xi}&=\lambda_-^\e x_- +\frac{\omega c_-^\e 2^\ell(1+k_c^2)}{\mu_c\kappa_0}(x_++x_-)|x_++x_-|^2 + \rmO_\e\left( |x_++x_-|^4\right)  .
\end{align}
\label{eqReduced3}
\end{subequations}

\subsection{Identification with the amplitude equation \eqref{eqAmpMod}}\label{sec:IDamplitude}

In this section, we will perform some linear transformations to place system \eqref{eqReduced3} into its normal form, which will, to leading order, be similar the amplitude equation \eqref{eqAmpMod} derived in the previous section. Using new variables $Y=x_++x_-$ and $Z=x_+-x_-$, we obtain
\begin{align*}
\dot Y &= -\frac{\alpha_1\epsilon}{2\alpha_2}Y+\frac{\sqrt{\Delta}\epsilon}{2\alpha_2}Z+\rmO\left(\e^2|Y+Z|+|Y+Z|^3+\e^{-1}|Y+Z|^4 \right), \\
\dot Z &= \frac{\sqrt{\Delta}\epsilon}{2\alpha_2}Y-\frac{\alpha_1\epsilon}{2\alpha_2}Z+\frac{2\omega(1+k_c^2)\alpha_2}{\mu_c\e  \sqrt{\Delta} \kappa_0}Y|Y|^2+\rmO\left(\e^2|Y+Z|+|Y+Z|^3+\e^{-1}|Y+Z|^4 \right).
\end{align*}
Rescaling time in the above system with $Y(t)=\epsilon u (\epsilon t)$ and $Z(t)=\epsilon v (\epsilon t)$ yields
\begin{align*}
\dot u &= -\frac{\alpha_1}{2\alpha_2}u+\frac{\sqrt{\Delta}}{2\alpha_2}v+\rmO(\e), \\
\dot v &= \frac{\sqrt{\Delta}}{2\alpha_2}u-\frac{\alpha_1}{2\alpha_2}v+\frac{2\omega(1+k_c^2)\alpha_2}{\mu_c \sqrt{\Delta}\kappa_0}u|u|^2+\rmO(\e).
\end{align*}
Finally, we make the transformation
\[u = q, \quad v = \frac{\alpha_1}{\sqrt{\Delta}}q + \frac{2\alpha_2}{\sqrt{\Delta}}p .\]
Then the equations take the form
\begin{subequations}
\begin{align}
\dot q &= p+\rmO(\e), \\
\dot p &= \frac{1}{\alpha_2}\left(-\alpha_0 q -\alpha_1 p +\frac{\omega(1+k_c^2)\alpha_2}{\mu_c\kappa_0}q|q|^2 \right)+\rmO(\e).
\end{align}
\label{systempq}
\end{subequations}
To leading order in $\e$, this system is equivalent to
\[ \ddot q + \frac{\alpha_1}{\alpha_2} \dot q +\frac{\alpha_0}{\alpha_2}q-\frac{\omega(1+k_c^2)}{\mu_c\kappa_0}q|q|^2=0.\]

Identifying these coefficients with the ones appearing in \eqref{eqAmpMod}, one can check via direct computations that the following equalities are satisfied
%we need to have:
\begin{align*}
-\frac{\hat\phi(k_c)}{\zeta}&=\frac{\alpha_0}{\alpha_2},\\
\frac{s}{\zeta}&=\frac{\alpha_1}{\alpha_2},\\
\frac{\omega}{\zeta}&=-\frac{\omega(1+k_c^2)}{\mu_c\kappa_0}.
\end{align*}
%This satisfied, as via direct computations, one gets
%\[ \zeta = \frac{4k_c^2(1+a^2+3k_c^2)}{(1+k_c^2)(a^2+k_c^2)}. \]
As a conclusion, system \eqref{systempq} is equivalent to \eqref{eqAmpModTF} that we rewrite as
\begin{subequations}
\begin{align}
\dot q &= p+\rmO(\e), \\
\dot p &= \frac{1}{\zeta}\left(\hat\phi(k_c) q -s p -\omega q|q|^2 \right)+\rmO(\e).
\end{align}
\label{systempqnew}
\end{subequations}

\subsection{Existence of heteroclinic orbits -- Proof of Theorem \ref{thmModTF}}\label{sec:heteroclinic}
We have transformed the reduced equation on the center manifold (\ref{eqReduced3}) into system \eqref{systempqnew}.  
To conclude the proof of Theorem~\ref{thmModTF} we will analyze this system and prove the existence of heteroclinic orbits corresponding to modulated traveling fronts.  To summarize, we show the existence of a circle of fixed points on the center manifold corresponding to the periodic solutions from Theorem~\ref{thmsolper} with differing phases.  When $\e=0$  the equation in the center manifold reduces to the real Ginzburg-Landau equation and heteroclinic connections between these fixed points and the fixed point at the origin are readily identified.  Then we show that these heteroclinic connections persist for $\e\neq 0$ and small.

We will refer to system \eqref{systempqnew} as $\mathcal{S}_\e$.  Then system $\mathcal{S}_0$  is
\begin{subequations}
\begin{align}
\dot q &= p, \\
\dot p &= \frac{1}{\zeta}\left(\hat\phi(k_c) q -s p -\omega q|q|^2 \right).
\end{align}
\label{systempqnew0}
\end{subequations}
We find that $\mathcal{S}_0$ has a fixed point at the origin  together with a circle of fixed points given by $p=0$, $|q|=\sqrt{\frac{\hat\phi(k_c)}{\omega}}$. Remember that $\hat\phi(k_c)<0$ and $\omega<0$ so that $\hat\phi(k_c)\omega>0$.  We have the following.
\begin{lemma}
For $s>\sqrt{-4\hat\phi(k_c)\zeta}$, the following assertions are satisfied.
\begin{enumerate}
\item[(i)] The origin $(q,p)=(0,0)$ is hyperbolic for $\mathcal{S}_\e$.
\item[(ii)] The system $\mathcal{S}_\e$ has a circle of normally hyperbolic fixed points which approach $|q|=\sqrt{\frac{\hat\phi(k_c)}{\omega}}$, $p=0$ as $\e\longrightarrow 0$.
\item[(iii)] For $\mathcal{S}_0$, for every $q$ on the circle $|q|=\sqrt{\frac{\hat\phi(k_c)}{\omega}}$, there is a saddle connection $\mathbf{C}_0$, tangent to the unstable direction at that point, which connects it to the origin $(q,p)=(0,0)$.
\item[(iv)] System $\mathcal{S}_\e$ has a family of heteroclinic connections $\mathbf{C}_\e$ (related to one another via $q\rightarrow e^{i\theta}q$ and $p\rightarrow e^{i\theta}p$) between the circle of fixed points and the origin.
\end{enumerate}
\end{lemma}
\begin{Proof}
\begin{itemize}
\item For (i), expand $p$ and $q$ in \eqref{systempqnew0} into their real and imaginary parts.  One finds that the linearization of the vector field at the origin $(q,p)=(0,0)$ has two double eigenvalues $\frac{-s\pm \sqrt{s^2+4\hat\phi(k_c)\zeta}}{2\zeta}$, which are both real and negative provided that $s>\sqrt{-4\hat\phi(k_c)\zeta}$ as $\zeta>0$. Thus, when $\e$ is small, the origin is also hyperbolic for $\mathcal{S}_\e$.
\item For (ii), consider first $\e=0$.  The linearization at $(p,q)=(0,\sqrt{\frac{\hat\phi(k_c)}{\omega}})$ has eigenvalues $0$, $\frac{-s}{\zeta}$ and $\frac{-s\pm \sqrt{s^2-8\hat\phi(k_c)\zeta}}{2\zeta}$. We thus have one unstable direction, one neutral direction and two stable directions. As a conclusion, for $\mathcal{S}_0$ the circle of fixed points $|q|=\sqrt{\frac{\hat\phi(k_c)}{\omega}}$, $p=0$ is normally hyperbolic.  When $\e$ is non-zero and small, the existence of small amplitude periodic solutions in Theorem~\ref{thmsolper} implies that this circle of fixed points persists as a circle of fixed points for $\mathcal{S}_\e$ which converges to $|q|=\sqrt{\frac{\hat\phi(k_c)}{\omega}}$ as $\e\to 0$.  This implies (ii).  
\item Statement (iii) is a well-known property of the real Ginzburg-Landau equation, see \cite{Aro-Wei}.  
\item Statement (iv) says that the heteroclinic solutions found in (iii) persist when $\e\neq 0$.  The proof of this fact is given in Lemma 4.2 of \cite{Eck-Way} and we do not repeat it here.  It relies on the normal hyperbolicity of the circle of fixed points and the relative dimensions of its unstable manifold and the stable manifold at the origin.  
\end{itemize}
\end{Proof}
We have thus shown that given $s>\sqrt{-4\hat\phi(k_c)\zeta}$, there is an $\e_0>0$ such that for all $\e\in (0,\e_0)$, \eqref{eq:NLKPP} has modulated traveling front solutions of frequency $k_c$ and of the form 
\[ u(t,x) = u(x-\e st,x)=\sum_{n\in \Z} W_n^u(x-\e st)e^{-ink_cx} \]
with the boundary conditions at infinity
\[ \underset{\xi\rightarrow - \infty}{\lim}u(\xi,x) = \mathbf{u}_\e(x)\approx 1+\e \sqrt{\frac{\hat\phi(k_c)}{\omega}}\cos(k_cx), \quad \underset{\xi\rightarrow + \infty}{\lim}u(\xi,x) =1.\]

\begin{rmk}
All our considerations above still apply if we replace $k_c$ with any $k$ for which $k=k_c+\delta$ and $\delta^2< -\frac{\hat\phi(k_c)}{1+\frac{\mu_c}{2}\hat\phi''(k_c)}\e^2$.
\end{rmk}

This concludes the proof of Theorem \ref{thmModTF}.

\section{Discussion}\label{sec:discussion}

We summarize our results and comment on extensions and major open questions.

\paragraph{Summary of results.} We studied the existence of stationary periodic solutions and modulated traveling fronts for nonlocal Fisher-KPP equations. Technically, we showed how center manifold reductions for infinite dimensional dynamical systems can elucidate the existence of such solutions for generic kernels in the case of stationary periodic solutions and for kernels with rational Fourier transform in the case of modulated traveling fronts. We also studied the spectral stability of the bifurcating stationary periodic solutions with respect to almost co-periodic perturbations using Bloch-wave decomposition.

\paragraph{Beyond exponential kernels.} Our results generalize, at least conceptually, to kernels with general rational Fourier transform. In those cases, one can still transform the nonlocal equation \eqref{eq:NLKPP} into a high-order system of partial differential equations of the form of \eqref{systemPDE}. One can still expect similar splitting phenomena for the eigenvalues of the resulting linear matrices $\mathcal{M}_n^\e$ so that a center manifold reduction along the lines of Eckmann \& Wayne \cite{Eck-Way} would apply.  

Extending our results to more general kernels would be more challenging.  On one hand, amplitude equations can be deduced by transforming (\ref{eq:NLKPP}) to Fourier space,
\[ \hv_t=\left(-k^2-\mu\hp(k)\right)\hv-\mu\hv\ast \hp\hv,\]
and then postulating a solution expansion of the form,
\[ \hv(t,k)=\sum_{|m|=1}^\infty 
\e^{|m|}\Phi_m(k,T)+\e^2\Phi_0(k,T).\]
This is referred to as a clustered mode distribution where the  $\Phi_j$ are $\rmO(1/\e)$ functions whose support lies on an $\rmO(\e)$ neighborhood of $k=j$, see for example \cite{Diprima,Mie2}.  Omitting the details, one can deduce a leading order equation for $\Phi_1$,
\[ \frac{\partial\Phi_1}{\partial T}=-\left(1+\frac{\mu_c}{2}\hat{\phi}''(k_c)\right)K^2\Phi_1-\hat{\phi}(k_c)\Phi_1+\omega \Phi_1\ast\Phi_1\ast\Phi_{-1},\]
where $K=(k-k_c)\e$.  This is the Fourier transformed version of the amplitude equation (\ref{eqAmpMod}).  Note that the same functional form holds in the general nonlocal case as for the case of exponential kernels studied here.  Let us also note that Morgan \& Dawes \cite{Mor-Daw} have successfully derived amplitude equations for the Swift-Hohenberg equation with nonlocal nonlinearity. 

The amplitude equation again suggests the existence of modulated traveling fronts in the genuinely nonlocal case.  It would be very valuable to adapt the center manifold techniques for modulated traveling waves to this context.  Some ideas on how to approach such questions in nonlocal problems can be found in \cite{Fay-Sch}.

\paragraph{Beyond small amplitude.} The numerical study of \cite{Nad-Per-Tan} and the results of Hamel \& Ryzhik \cite{Ham-Ryz} show that the periodic stationary solutions constructed in the paper persist beyond $\mu>\mu_c$, with amplitude that no longer scales as $\sqrt{\mu-\mu_c}$. It would be interesting to use general continuation techniques to study how these periodic stationary solutions perturb as we vary the parameter $\mu$ for example. For example, could such a branch of solutions reconnect somewhere to the solution $u=1$ is a question that we would like to address in future work.

\paragraph{Stability.} Another problem of interest is the question of stability of the periodic solutions found in this paper with respect to either co-periodic or general non-periodic perturbations. Our results in Theorem \ref{thmstability} only give a partial answer to the problem as we are only able to characterize regions of spectral instability with respect to almost co-periodic perturbations. Stability of the bifurcating periodic solutions has been well studied  in the context of the Swift-Hohenberg equation \cite{Mie,Mie1}.  There criteria have been derived namely which give the region in parameter space $(\e,\delta)$ where roll solutions of the Swift-Hohenberg equation are spectrally stable \cite{Mie}. The analysis uses a Bloch-wave decomposition to study the spectrum of linearized equation about a periodic solution.    It would then be interesting to adapt such techniques in the context of the  nonlocal problem (\ref{eq:NLKPP}).

\paragraph{Two-stage invasion fronts.} From the perspective of the original problem (\ref{eq:NLKPP}) and the related Fisher-KPP equation it is often the dynamics for initial data near the state $u=0$ that is of interest.  Here, one observes traveling fronts where the zero state is invaded by a periodic stationary state around $u=1$.  Sometimes an intermediate region where the solution is approximately in the state $u=1$ is observed.   Invasion fronts of this form were numerically computed in \cite{Nad-Per-Tan}.  Since the stationary periodic solutions come in families, one expects that the invasion process is dynamically selecting a particular pattern amongst this family of solutions.  When $\mu\approx\mu_c$, the primary front where $u=1$ replaces $u=0$ travels much faster than the secondary modulated front and the selected pattern is determined by the modulated traveling front propagating with the minimal speed.  However, when $\mu$ is large the numerically observed speeds of the secondary modulated traveling front exceed or are of the same order as that of the primary front and the pattern selection mechanism is more difficult to characterize.  This will be the object of future study.

\section*{Acknowledgments}
We would like to thank the referee for pointing out a gap in the application of the center manifold result of Proposition \ref{prop:cm}.  The research of GF leading to these results has received funding from the European Research Council under the European Union's Seventh Framework Program (FP/2007-2013) / ERC Grant Agreement n321186 : "Reaction-Diffusion Equations, Propagation and Modelling".

\appendix

\section{Computation of $\langle \psi_\pm^\e, \varphi_\pm^\e \rangle$}\label{AppIP}

To simplify notation, we denote by $\mathcal{M}_\e$ the matrix $\mathcal{M}_1^\e$. From its definition (see equation from \eqref{eqMatrixMn}), $\mathcal{M}_\e$ can be decomposed as
\[\mathcal{M}_\e= \mathcal{M}_0+\e \mathcal{L}_1 +\e^2 \mathcal{L}_2,  \]
where we have
\[
\mathcal{M}_0:=\left(\begin{matrix} 0 & 1 & 0&0 & 0 & 0 \\
k_c^2 & 2 i k_c & \mu_c & 0 & \mu_c & 0 \\
0 & 0 & 0&1 & 0 & 0 \\
-3a^2 & 0 & a^2+k_c^2 & 2 i  k_c & 0 & 0 \\
0 & 0 & 0&0 & 0 & 1 \\
2 & 0 & 0 & 0 & 1+k_c^2 & 2 i  k_c \end{matrix}\right),\]
and 
\[\mathcal{L}_1:=\left(\begin{matrix} 0 & 0 & 0&0 & 0 & 0 \\
0 &  -s & 0 & 0 & 0 & 0 \\
0 & 0 & 0&0 & 0 & 0 \\
0 & 0 & 0 & 0 & 0 & 0 \\
0 & 0 & 0&0 & 0 & 0 \\
0 & 0 & 0 & 0 & 0 & 0 \end{matrix}\right), \quad \mathcal{L}_2:=\left(\begin{matrix} 0 & 0 & 0&0 & 0 & 0 \\
0 &  0 & 1 & 0 & 1 & 0 \\
0 & 0 & 0&0 & 0 & 0 \\
0 & 0 & 0 & 0 & 0 & 0 \\
0 & 0 & 0&0 & 0 & 0 \\
0 & 0 & 0 & 0 & 0 & 0 \end{matrix}\right).\]
We define $\rme_0$ and $\rme_1$, the eigenvector and generalized eigenvector of $\mathcal{M}_0$ such that $\mathcal{M}_0\rme_0=0$ and $\mathcal{M}_0\rme_1=\rme_0$. Recall that
\[\rme_0:= \left(1,0,\hat\phi_v(k_c),0,\hat\phi_w(k_c),0\right)^T, \quad \rme_1:=\left( 0, 1, -\frac{2ik_c}{a^2+k_c^2}\hat\phi_v(k_c),\hat \phi_v(k_c),-\frac{2ik_c}{1+k_c^2}\hat\phi_w(k_c),\hat\phi_w(k_c)\right)^T. \]
Next, we define the adjoint matrix to $\mathcal{M}_\e$ as
\[ \mathcal{M}_\e^*= \mathcal{M}_0^*+\e \mathcal{L}_1 +\e^2 \mathcal{L}_2^T, \]
where
\[ \mathcal{M}_0^*=\left(\begin{matrix} 
0 & k_c^2 & 0 & -3a^2 & 0 & 2 \\
1 & -2ik_c & 0 & 0 & 0 & 0 \\
0 & \mu_c & 0 & a^2+k_c^2 & 0 & 0 \\
0 & 0 & 1 & -2ik_c & 0 & 0 \\
0 & \mu_c & 0 & 0 & 0 & 1+k_c^2 \\
0 & 0 & 0 & 0 & 1 & -2ik_c
\end{matrix} \right). \]
We define $\rme_0^*$ and $\rme_1^*$, the eigenvector and generalized eigenvector of $\mathcal{M}_0^*$ such that $\mathcal{M}_0^*\rme_0^*=0$ and $\mathcal{M}_0^*\rme_1^*=\rme_0^*$. Furthermore, we impose the orthogonality conditions
\[ \langle \rme_0 , \rme_0^*\rangle =0, \quad \langle \rme_0 , \rme_1^*\rangle =1,\quad \langle \rme_1 , \rme_1^*\rangle =0, \quad \langle \rme_1 , \rme_0^*\rangle =1. \]
We have
\begin{align*} 
\rme_0^*&:=\frac{1}{\kappa_0}\left(-\frac{2ik_c(1+k_c^2)}{\mu_c},-\frac{1+k_c^2}{\mu_c},\frac{2ik_c(1+k_c^2)}{a^2+k_c^2},\frac{(1+k_c^2)}{a^2+k_c^2}, 2ik_c,1\right)^T,\\
\rme_1^*&:=\frac{1}{\kappa_1}\left(-\frac{1+5k_c^2}{\mu_c},\frac{2ik_c}{\mu_c},\frac{a^2+5a^2k_c^2-3k_c^2+k_c^4}{(a^2+k_c^2)^2},\frac{2ik_c(1-a^2)}{(a^2+k_c^2)^2},1,0\right)^T,
\end{align*}
where we have set
\begin{eqnarray}
\kappa_0&:=&-4k_c^2(1+k_c^2)\left(\frac{3a^2}{(a^2+k_c^2)^3}-\frac{2}{(1+k_c^2)^3} \right)=-\zeta\frac{(1+k_c^2)}{\mu_c}, \label{eq:kappa0} \\
\kappa_1&:=& -\frac{1+5k_c^2}{\mu_c}+\hat\phi_v(k_c)\frac{a^2+5a^2k_c^2-3k_c^2+k_c^4}{(a^2+k_c^2)^2}+\hat\phi_w(k_c). \nonumber
\end{eqnarray}

We recall that we defined $\varphi_\pm^\e$ as the corresponding eigenvectors to $\mathcal{M}_\e$ that satisfy
\[ \mathcal{M}_\e \varphi_\pm^\e = \lambda_\pm^\e \varphi_\pm^\e, \quad \lambda_\pm^\e=\e\chi_\pm +\rmO(\e^2). \]
Similarly, we have $\psi_\pm^\e$ the eigenvectors of $\mathcal{M}_\e^*$ that satisfy
\[ \mathcal{M}_\e^* \psi_\pm^\e = \bar \lambda_\pm^\e \psi_\pm^\e, \quad \bar \lambda_\pm^\e=\e\chi_\pm +\rmO(\e^2). \]

We want to evaluate, at leading order in $\e$, the scalar product $\langle \psi_\pm^\e, \varphi_\pm^\e \rangle$. To do so, we write an expansion for $\varphi_\pm^\e$ and $\psi_\pm^\e$ of the form
\[ \varphi_\pm^\e = \rme_0 +\e \xi_\pm + \rmO(\e^2) \text{ and } \psi_\pm^\e = \rme_0^* +\e \xi_\pm^* + \rmO(\e^2). \]
And then, we have
\[ \langle \psi_\pm^\e, \varphi_\pm^\e \rangle = \left( \langle \rme_0^*,\xi_\pm \rangle + \langle \xi_\pm^* , \rme_0 \rangle \right)\e+\rmO(\e^2). \]
One needs to compute the two inner products $\langle \rme_0^*,\xi_\pm \rangle$ and $\langle \xi_\pm^* , \rme_0 \rangle$. For this, we expand the relation $\mathcal{M}_\e \varphi_\pm^\e = \lambda_\pm^\e \varphi_\pm^\e$ and obtain
\[ \left( \mathcal{M}_0+\e \mathcal{L}_1 +\e^2 \mathcal{L}_2 \right) \left( \rme_0 +\e \xi_\pm + \rmO(\e^2) \right) = (\e\chi_\pm + \rmO(\e^2))(\rme_0 +\e \xi_\pm + \rmO(\e^2)). \]
Collecting terms of order $\rmO(1)$, we have $\mathcal{M}_0\rme_0=0$ which is satisfied by definition of $\rme_0$. At order $\rmO(\e)$, we obtain an equation
\[ \mathcal{M}_0 \xi_\pm + \mathcal{L}_1\rme_0 = \chi_\pm \rme_0, \]
and taking the inner product with $\rme_1^*$ yields
\[ \langle \xi_\pm , \rme_0^*\rangle = \chi_\pm \langle \rme_0 , \rme_1^* \rangle - \langle \mathcal{L}_1\rme_0 , \rme_1^*\rangle=\chi_\pm, \]
as $\mathcal{L}_1\rme_0=0$. Similar computations for $\psi_\pm^\e$ give
\[ \mathcal{M}_0^* \xi_\pm^* + \mathcal{L}_1\rme_0^* = \chi_\pm \rme_0^*, \]
and taking the inner product with $\rme_1$, we obtain
\[ \langle \xi_\pm^* , \rme_0 \rangle = \chi_\pm \langle \rme_0^* , \rme_1 \rangle - \langle \mathcal{L}_1\rme_0^* , \rme_1\rangle=\chi_\pm -\frac{(1+k_c^2)s}{\mu_c\kappa_0}. \]
As a conclusion, we have that
\[ \langle \psi_\pm^\e, \varphi_\pm^\e \rangle = \left(  2\chi_\pm -\frac{(1+k_c^2)s}{\mu_c\kappa_0}   \right)\e+\rmO(\e^2).\]
Finally, we can check that
\[ -\frac{\alpha_1}{\alpha_2} -\frac{(1+k_c^2)s}{\mu_c\kappa_0} =0,\]
 as $\frac{s}{\zeta}=\frac{\alpha_1}{\alpha_2}$ and $\kappa_0=-\zeta\frac{1+k_c^2}{\mu_c}$, such that
\be \langle \psi_\pm^\e, \varphi_\pm^\e \rangle = \pm \frac{\sqrt{\Delta}}{\alpha_2}\e+\rmO(\e^2).
\label{eqInnerProduct}
\ee

\bibliography{plain}

\end{document}